\RequirePackage{fix-cm}
\documentclass[smallextended]{svjour3}       
\smartqed  
\usepackage{graphicx}
\usepackage{amsfonts, amssymb}
\usepackage{amsmath}
\usepackage{graphicx}
\usepackage{tikz}
\usepackage{tikzscale}
\usepackage{xcolor}
\usepackage{color}
\usepackage{xspace}
\usepackage{calrsfs}
\usepackage{algorithm}
\usepackage{algpseudocode}
\usepackage{adjustbox}
\usepackage{xspace}
\usepackage{sidecap}
\usepackage[bookmarks,pagebackref,
pdfpagelabels=true, 
]{hyperref}

\usepackage{makeidx}
\usepackage{imakeidx}
\makeindex[columns=1]
\newcommand{\LieG}{\mathcal{G}}
\newcommand{\LieGhat}{\hat{\mathcal{G}}}
\newcommand{\LieAlg}{\mathcal{L}}
\newcommand{\LieAlghat}{\hat{\mathcal{L}}}
\newcommand{\Sys}{R}
\newcommand{\Syshat}{\hat R}

\newcommand{\dd}[1]{\frac{\partial}{\partial #1}}

\newcommand{\BK}{\sc{BK}}
\newcommand{\ID}{\sc{ID}}

\newcommand{\DE}{DE\xspace}
\newcommand{\DEs}{DEs\xspace}
\newcommand{\PDE}{PDE\xspace}
\newcommand{\ODE}{ODE\xspace}

\newcommand{\DPS}{DPS\xspace}

\newcommand{\xh}{\hat{x}}
\newcommand{\vh}{\hat{v}}
\newcommand{\uh}{\hat{u}}
\newcommand{\thn}{\hat{t}}

\newcommand{\alg}[1]{{\sc #1}\xspace}
\newcommand{\maple}[1]{{\tt #1}\xspace}
\newcommand{\inputs}[1]{{\it #1}\xspace}
\newcommand{\MapDE}{\maple{MapDE}}
\newcommand{\PreEquivTest}{\maple{PreEquivTest}}

\newcommand{\Source}{\inputs{Source}}

\newcommand{\Target}{\inputs{Target}}

\newcommand{\Map}{\inputs{Map}}

\newcommand{\Maple}{\maple{Maple}}
\newcommand{\LAVF}{\maple{LAVF}}

\newcommand{\rif}{\alg{rif}}
\newcommand{\dec}{\alg{dec}}

\newboolean{showComments}
\newcommand{\todocmd}[1]{\small{\textcolor{red}{#1}}}

\newcommand{\todo}[1]{\ifthenelse {\boolean{showComments}} {\todocmd{#1}} {}}
\begin{document}

\title{Symmetry-based algorithms for invertible mappings of polynomially nonlinear PDE to linear PDE {\todo{[R1C1]}}}
%


\titlerunning{Symmetry-based algorithms for mapping nonlinear PDE to linear PDE}        

\author{Zahra Mohammadi \and
        Gregory J. Reid \and S.-L. Tracy Huang
}


\institute{Zahra Mohammadi, Gregory J. Reid \at
              Department of Applied Mathematics, University of Western Ontario, Canada \\
              \email{zmohamm5@uwo.ca; reid@uwo.ca}           
\and
S.-L. Tracy Huang \at
Data61, CSIRO, Canberra ACT 2601, Australia \\
\email{ Tracy.Huang@data61.csiro.au  }
}

\date{Received: date / Accepted: date}

\maketitle

\begin{abstract}
This paper is a sequel to our previous work where we introduced the \MapDE algorithm to determine the existence of analytic invertible mappings of an input (source) \index{{\DPS}, Differential Polynomial System} differential polynomial system ({\DPS}) to a specific target \DPS,
and sometimes by heuristic integration an explicit form of the mapping.
A particular feature was to exploit the Lie symmetry invariance algebra of the source without integrating its equations, to facilitate 
{\MapDE}, making algorithmic an approach initiated by Bluman and Kumei.
In applications, however, the explicit form of a target {\DPS} is not available, and a more important question is, can the 
source be mapped to a more tractable class?  

We extend {\MapDE} to determine if a source nonlinear {\DPS} can be mapped to a linear differential system.
\MapDE applies differential-elimination completion algorithms to the various over-determined {\DPS} by applying a finite number of differentiations and eliminations to complete them to a form for which an existence-uniqueness theorem is available, enabling the existence of the linearization to be determined among other applications.
The methods combine aspects of the Bluman-Kumei mapping approach with techniques introduced by Lyakhov, Gerdt and Michels for the determination of exact linearizations of \index{{\ODE}, Ordinary Differential Equation} {\ODE}.
The Bluman-Kumei approach for  \index{{\PDE}, Partial Differential Equation} {\PDE} focuses on the fact that such linearizable systems must admit
a usually infinite Lie sub-pseudogroup corresponding to the linear superposition of solutions in the target.
In contrast, Lyakhov et al. focus on {\ODE} and properties of the so-called derived sub-algebra of the (finite) dimensional Lie algebra of symmetries of the {\ODE}.  
Examples are given to illustrate the approach, and a heuristic integration method sometimes gives explicit forms of the maps.
We also illustrate the powerful maximal symmetry groups facility as a natural tool to be used in conjunction with {\MapDE}. 



\end{abstract}

\section{Introduction}
\label{sec:Intro}
This paper is a sequel to \cite{MohReiHua19:Intro} and is part of a series in which we explore algorithmic aspects of exact and approximate mappings of differential equations.  We are interested in mapping less tractable differential equations into more tractable ones, in particular in this article focusing on mapping nonlinear systems to linear systems.  This builds on progress in
\cite{MohReiHua19:Intro} where we considered mappings from a specific differential system to a specific target system and mappings from a linear to a linear constant coefficient differential equation.

As in \cite{MohReiHua19:Intro} we consider systems of (partial or ordinary) differential equations with $n$ independent variables and 
$m$ dependent variables which are local analytic functions of their arguments. Suppose $\Sys$ has independent variables $x = (x^1, \ldots ,x^n )$ and dependent variables $u = (u^1, \ldots , u^m)$ and 
$\Syshat$ has independent variables $\hat{x} = (\hat{x}^1, \ldots , \hat{x}^n )$ and dependent variables $\hat{u} = (\hat{u}^1, \ldots , \hat{u}^m)$.
In particular, we consider local analytic mappings $\Psi$:
$(\hat{x}, \hat{u} ) = \Psi (x,u) = (\psi (x,u), \phi(x,u) )$ so that $\Sys$ is locally and invertibly mapped to $\Syshat$:
\begin{equation}
	\label{eq:MapTrans}
	\hat{x}^j  =   \psi^j (x,u), \qquad 
	\hat{u}^k =   \phi^k (x,u) 
\end{equation}
where $j = 1,\ldots, n$ and $k = 1,\ldots, m$.
The mapping is locally invertible, so \index{$\mbox{Det} \mbox{Jac}$, Determinant of the Jacobian} the determinant of the Jacobian of the mapping is nonzero: 
\begin{equation}
	\label{eq:Jac}
	\mbox{Det} \mbox{Jac}(\Psi) = \mbox{Det} \frac{\partial (\psi, \phi )}{\partial (x, u)} \not = 0,
\end{equation}
where $\frac{\partial (\psi, \phi )}{\partial (x, u)} $ is the usual Jacobian $(n+m) \times (n+m)$ matrix of first order derivatives of the $(n+m)$ functions $(\psi, \phi)$ with respect to the $(n+m)$ variables $(x, u)$.
Note throughout this paper: we will call $\Syshat$ the 
\Target system of the mapping, which will generally have some more desirable features than $\Sys$, which we call the \Source system.
A well-known direct approach to forming the equations satisfied by $\Psi$ is roughly to substitute the general change of variables into $\Syshat$, evaluate the result modulo $\Sys$,
appending equations that express the independence of $\Psi$ on derivative jet variables
(or equivalently decomposing in independent expressions in the jet variables).
The resulting equations for $\Psi$ are generally nonlinear overdetermined systems.
Algorithmic manipulation of these and other over-determined systems of {\PDE} are at the core of the algorithms
used in this paper.  We make prolific use of \index{{\dec}, Differential-elimination completion algorithms} differential-elimination completion  ({\dec}) algorithms, {\todo{[R1C3]}} which 
apply a finite number of differentiations and eliminations to complete such over-determined systems to a form including their integrability conditions, for which an existence uniqueness theorem is available.  {\Maple} is fortunate to have several such differential elimination packages. 
Currently we use  the {\rif} algorithm via the {\Maple} command $\maple{rifsimp}$ \cite{Rus99:Exi} in our implementation, but other \Maple packages could be used \cite{Robertz106:Tdec,BLOP:Diffalg}.
 {\todo{[R1C2]}} To be algorithmic we restrict to systems $R$ and $\hat{R}$ that are polynomially nonlinear (i.e. differential polynomial systems, DPS).
In this  paper \dec refers to a Differential Elimination Completion algorithm
to emphasize that a number of algorithms are available.

A very general approach to such problems, concerning maps $\Psi$ from $\Sys$ to $\Syshat$, is Cartan's famous Method of Equivalence which
finds invariants that label the classes of systems equivalent under the pseudogroup of such mappings. See especially texts \cite{PetOlv107:Sym} and  \cite{Man10:Pra}. The fundamental importance and computational difficulty of such equivalence questions 
has attracted attention from the symbolic computation community \cite{NeutPetitotDridi2009}.
For recent developments and extensions of Cartan's moving frames for equivalence problems see \cite{Fel99:Mov}, \cite{Valiquette13:LocEquiv} and \cite{Arnaldsson17:InvolMovingFrames}.
The \maple{DifferentialGeometry} package \cite{And12:New} is available in {\Maple} and has been applied to equivalence problems \cite{KruglikovThe18}.
Underlying these calculations is that overdetermined \PDE systems with some 
non-linearity must be reduced to forms that enable the statement of a local existence and uniqueness theorem \cite{Hub09:Dif,GolKonOvcSza09:complexityDiffElim,Sei10:Inv,Rus99:Exi,Bou95:Rep}.

Our methods here and in \cite{MohReiHua19:Intro} are based on the mapping approach initiated by Bluman and Kumei~\cite{BluKu109:DiEq} which focuses on the interaction between such mappings
and Lie symmetries via their infinitesimal form on the source and target.
In particular, let $\LieG$ be the Lie group of transformations leaving $\Sys$ invariant.{\todo{[R4C2]}}
Also, let $\LieGhat$ be the Lie group of transformations leaving $\Syshat$ invariant.
Locally, such Lie groups are characterized by their linearizations in a neighborhood of their
identity, that is by their Lie algebras $\LieAlg$, $\LieAlghat$.
If an invertible map $\Psi$ exists then $\LieG \simeq \LieGhat$ and $\LieAlg \simeq \LieAlghat$.
This yields a subsystem of linear equations for $\Psi$ which we call the Bluman-Kumei equations.
{\todo{[R4C3]}}
It is a significant challenge to translate the methods of Bluman and 
Kumei into procedures that are algorithmic (i.e.,\ guaranteed to succeed on a defined class of inputs in finitely many steps).
Please see ~\cite{AnBlWo110:Mapp,BLuYang119:SysAlg,Wolf116:SymSof,TWolf108:ConLa}
for progress in their approach and some (heuristic) integration-based computer implemented methods.
Our methods are also inspired by remarkable recent progress on this question for \ODE by Lyakhov et al.\cite{LGM101:LG} who presented an algorithm for determining when an \ODE is linearizable.
It was also stimulated by their use of an early method by one of us (see Reid \cite{Rei91:Fin}), which has been dramatically improved and extended \cite{RLB92:Alg,LisleReidInfinite98}
with the latest improvements in the {\LAVF} package \cite{LisH:Alg,Huang:Thesis}.

In our previous work \cite{MohReiHua19:Intro} we provided an algorithm to determine the existence of a mapping of a linear differential equation to the class of constant coefficient linear homogeneous differential equations.  Key for this application was the exploitation of a commutative
sub-algebra of symmetries of $\LieAlghat$ corresponding to translations of the independent variables in the target. 
The main contribution of this paper is to present an algorithmic method for determining the mapping of a nonlinear system to a linear system when it exists.  
Using a technique of Bluman and Kumei, we exploit the fact that $\Syshat$ must admit a sub-pseudo group corresponding to the superposition property that linear systems by definition must satisfy.  Once existence is established, a second stage can determine features of the map and sometimes by integration, explicit forms of the mapping.
For an algorithmic treatment using differential elimination (differential algebra), we limit our treatment to  systems of differential polynomials, with coefficients from $\mathbb{Q}$ or some computable extension of $\mathbb{Q}$ in $\mathbb{C}$.  Thus our input system $\Sys$ should be a system of {\DPS}.
 Some non-polynomial systems can be converted to differential polynomial form by using the Maple command, \maple{dpolyform}.

In \S \ref{sec:GeomDPSID} we provide some introductory material on differential-elimination algorithms, initial data and Hilbert dimensions. In \S \ref{sec:PreMapEqs} we give an introduction to symmetries and mapping equations.  In \S \ref{sec:Algorithms}, we 
introduce the \MapDE algorithm.  Examples of application \MapDE are given in \S \ref{sec:Examples}, and we conclude with  a discussion in \S \ref{sec:Discussion}. Our \MapDE program and a demo file are publicly available on GitHub at: 
 \href{https://github.com/GregGitHub57/MapDETools}{https://github.com/GregGitHub57/MapDETools}.

\section{Differential-elimination algorithms, initial data and Hilbert functions}

\label{sec:GeomDPSID}

The geometric approach to {\DPS} centers on the jet locus, the solution set of the equations obtained by replacing derivatives with formal variables, yielding systems of polynomial
equations and inequations and differences of varieties (solution sets of polynomial equations).
In this way the algorithmic tools of algebraic geometry can be applied to systems of {\DPS}.
The union of prolonged graphs of local solutions of a {\DPS} is a subset of the jet locus in $J(\mathbb{C}^n, \mathbb{C}^m)$, the jet space, with $n$ independent variables, and $m$ dependent variables.
For details concerning Jet geometry see \cite{Olv93:App,Sei10:Inv}. 

Throughout this paper we make prolific use of differential-elimination algorithms which apply a finite number
of differentiations and eliminations to an input {\DPS} to yield in a form that yields information about its 
properties and solutions.  For example, consider 
\begin{equation}
\label{EzSys}
u_{xyy} - u_{yy} = 0,   \; \; u_{xx} + u_{xy} - u_x - u_y = 0,  \; \;  u_{xx} - u_{xy} - u_x + u_y = 0
\end{equation}
Simply eliminating using the ordering 
$u_{xx} \succ u_{xy} \succ u_{yy} \succ u_{x}  \succ u_{y}$
gives the equivalent system $u_{xyy} = u_{yy} ,   u_{xx} = u_x ,   u_{xy} = u_y$.
The first 
equation can be omitted since it is a derivative of the third, yielding 
\begin{equation}
\label{EzSys2}
u_{xx} = u_x , \; \; \;   u_{xy} = u_y
\end{equation}
The operations to reduce the example above mirror those to reduce a related polynomial system via
$\frac{\partial}{\partial x} \leftrightarrow X$, $\frac{\partial}{\partial y} \leftrightarrow Y$, 
$X Y^2 - Y^2 = 0$, $X^2  + X Y - X - Y = 0$, $X^2  - X Y - X  + Y = 0$ to a Gr\"obner Basis.
Indeed a natural generalization of Gr\"obner Bases exists for linear homogeneous {\PDE}.
However {\DPS} are much tougher theoretically and computationally, with straightforward 
generalizations yielding infinite bases, and undecidable problems.
Currently we use  the {\rif} algorithm via the {\Maple} command $\maple{rifsimp}$ \cite{Rus99:Exi} in our implementation, but other \Maple packages could be used
 such as \maple{DifferentialThomas} Package \cite{Thomas106:Tdec,Robertz106:Tdec},{\todo{[R1C4]}} and the \maple{DifferentialAlgebra} package\cite{Bou95:Rep,BLOP:Diffalg,lemaire:tel} {\todo{[R1C5]}} or \maple{casesplit} which offers a uniform interface to such packages.

A key aspect of these {\dec} packages for {\DPS} is that they split on cases where certain leading polynomial quantities are zero or nonzero.  This leads to systems of differential polynomial equations and inequations.
In particular a system $R$ of equations $\{ p_1 = 0, p_2 = 0, \cdots , p_b =0 \} $
and inequations  $ \{ q_1 \not = 0, \cdots , q_c \not = 0 \}$
has solution locus
\begin{equation}
\label{Z=E-I}
Z(R) = V^= (R) \setminus V^{\not =} (R)
\end{equation}
where 
$V^= (R)$ are the solutions satisfying $\{ p_1 = 0, p_2 = 0, \cdots , p_c =0 \} $
and  $V^{\not =}(R)$ is the set of solutions of $\Pi_{i = 1}^{i = c} \: q_i = 0$.

Moreover, a central input in such algorithms are rankings of derivatives \cite{Rus97:Ran}.  Indeed 
let $\Omega(\Sys)$ be all the derivatives of dependent variables for $\Sys$.
Throughout this paper the set of derivatives also includes $0$-order derivatives (i.e. dependent variables).
A ranking on $\Omega(\Sys)$ is a total order $\prec$ that satisfies the axioms in \cite{Rus97:Ran}.
Given a ranking and algorithms in \cite{Rus99:Exi} determine initial data and the existence and uniqueness of formal
power series solutions. Additionally, if the ranking is orderly and of Riquier type (i.e. ordered first by total order of derivative, with a ranking specified by a Riquier ranking matrix) analytic initial data yields local analytic solutions. See \cite{Riquier:diffalg} {\todo{[R1C6]}} for a proof of this result.  We will need some block elimination rankings that eliminate groups of dependent variables in favor of others. Enforcing the block order via the first row of the Riquier Matrix and then enforcing total order of the derivative as the next criterion for each block enables analytic data to yield analytic solutions that is sufficient for this 
paper.

The differential-elimination algorithms used in this article enable the algorithmic posing of initial data for the
determination of unique formal power series of differential systems.  We will exploit a powerful measure of solution dimension information given by Differential Hilbert Series and its related 
Differential Hilbert Function \cite{MLPK13:Hilbert,MLH:DCP}.
Indeed given a ranking algorithms such as the \rif algorithm, Rosenfeld Groebner and Thomas Decomposition partition the set of derivatives of the unknown function at a regular point $x_0$ into a set of parametric derivatives $\mathcal{P}$ \index{$\mathcal{P}$, A set of parametric derivatives} and a set of principal derivatives.  Principal derivatives are derivatives of the leading derivatives and 
parametric derivatives are its complement, the ones that can be ascribed arbitrary values at $x_0$.
Then the Hilbert Series is defined as:
\begin{equation}
\label{HS1}
\mbox{HS}(\Sys, x_0) := \sum_{G \in  \mathcal{P} } s^{\mbox{dord}(G)} = \sum_{n=0}^\infty a_n s^n
\end{equation}
\index{$\mbox{dord}$, Differential Order} where $\mbox{dord}(G)$ is the differential order of $G$, $a_n$ is the number of parametric derivatives of differential order $n$. 
To algorithmically compute such a series exploits the fact that the parametric data can be partitioned into 
subsets.

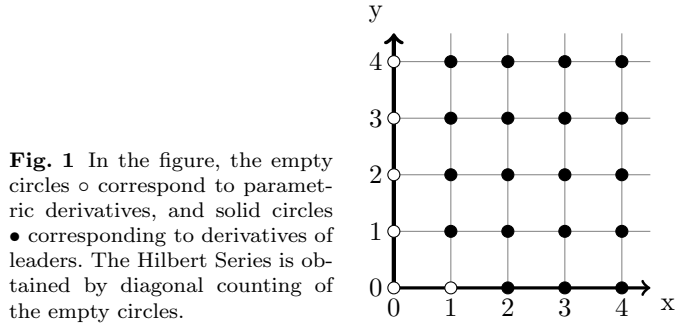
\begin{SCfigure}
	\label{fig:HilbertPic}
	\begin{tikzpicture}[scale=.75]
	\draw[step=1cm,gray,very thin] (0,0) grid (4.5,4.5);
	\draw[very thick,->] (0,0) -- (4.5,0) node[anchor=north west] {x};
	\draw[ultra thick,->] (0,0) -- (0,4.5) node[anchor=south east] {y};
	\foreach \x in {0,1,2,3,4}
	\draw (\x cm,1pt) -- (\x cm,-1pt) node[anchor=north] {$\x$};
	\foreach \y in {0,1,2,3,4}
	\draw (1pt,\y cm) -- (-1pt,\y cm) node[anchor=east] {$\y$};
	\draw [fill=black] (1,1) circle[radius= 0.3 em]; 
	\draw [fill=black] (1,2) circle[radius= 0.3 em];
	\draw [fill=black] (1,3) circle[radius= 0.3 em];
	\draw [fill=black] (1,4) circle[radius= 0.3 em];
	\draw [fill=black] (2,0) circle[radius= 0.3 em];
	\draw [fill=black] (2,1) circle[radius= 0.3 em];
	\draw [fill=black] (2,2) circle[radius= 0.3 em];
	\draw [fill=black] (2,3) circle[radius= 0.3 em];
	\draw [fill=black] (2,4) circle[radius= 0.3 em];
	\draw [fill=black] (3,0) circle[radius= 0.3 em];
	\draw [fill=black] (3,1) circle[radius= 0.3 em];
	\draw [fill=black] (3,2) circle[radius= 0.3 em];
	\draw [fill=black] (3,3) circle[radius= 0.3 em];
	\draw [fill=black] (3,4) circle[radius= 0.3 em];
	\draw [fill=black] (4,0) circle[radius= 0.3 em];
	\draw [fill=black] (4,1) circle[radius= 0.3 em];
	\draw [fill=black] (4,2) circle[radius= 0.3 em];
	\draw [fill=black] (4,3) circle[radius= 0.3 em];
	\draw [fill=black] (4,4) circle[radius= 0.3 em];
	\draw [fill=white] (0,0) circle[radius= 0.3 em]; 
	\draw [fill=white] (0,1) circle[radius= 0.3 em]; 
	\draw [fill=white] (0,2) circle[radius= 0.3 em]; 
	\draw [fill=white] (0,3) circle[radius= 0.3 em]; 
	\draw [fill=white] (0,4) circle[radius= 0.3 em]; 
	\draw [fill=white] (1,0) circle[radius= 0.3 em]; 
	\end{tikzpicture}
	\caption{In the figure, the empty circles $\circ$ correspond to parametric derivatives, and solid circles $\bullet $ corresponding to derivatives of leaders.  The Hilbert Series is obtained by diagonal counting of the empty circles.}
\end{SCfigure}

For example consider $\Sys = \{ u_{xx} = u_x, u_{xy} = u_y\}$, which is already in \rif-form with respect
to an orderly ranking.  Then $\mathcal{P} = \{ u_x \} \cup \{ u , u_y, u_{yy}, u_{yyy}, \cdots \}$, and the 
associated set of initial data is $ \{ u_x (x_0, y_0) = c_0 \} \cup \{ u (x_0, y_0) = c_1  , u_y (x_0, y_0) = c_2 , u_{yy} (x_0, y_0) = c_3 , \cdots \}$.
In what follows it is helpful to associate these derivatives with corresponding points in $\mathbb{N}^2$ via
$ \frac{\partial^{i+j}}{\partial x^i  \partial y^j} u \leftrightarrow (i, j) \in \mathbb{N}^2$.
See  Fig.~\ref{fig:HilbertPic} for a graphical depiction of $\mathcal{P}$.

Then 
\begin{equation}
\label{HS2}
\mbox{HS}(\Sys, (x_0, y_0)) =  1 + 2s + s^2 + s^3 + s^4 + \cdots
\end{equation}
A crucial condition to check in our algorithms will be that two Hilbert series are the same, which is
complicated since no finite algorithm exists for checking equality of series.
For our example, however, the series can be expressed finitely
$\mbox{HS}(\Sys, (x_0, y_0)) =  1 + 2s + s^2 + s^3 + s^4 + \cdots = s + \frac{1}{1 - s}$.
This collapsing of the series into a rational function can be accomplished in general.
In our implementation, we exploited the output of the \maple{initialdata} algorithm in the {\rif} package.
For example this returns data as a partition of two sets: a finite set of initial data $\mathcal{F}$ and an
set that represents an infinite set of data $\mathcal{I}$ \index{$\mathcal{I}$, An infinite set of data} by compressing them into arbitrary functions.
For our example this compression is:
\begin{equation}
\label{FuI}
 \mathcal{F} \cup \mathcal{I} = \{ u_x (x_0, y_0) = c_0 \} \cup \{ u (x_0, y) = f(y) \}
\end{equation}
This approach is easily extended to yield the Differential Hilbert Function \index{$\mbox{HF}$,Differential Hilbert Function} by applying a function that acts on each 
piece of initial data in $\mathcal{F}$ and $\mathcal{I}$
\begin{equation}
\label{HF(FuI)}
\mbox{HF} (\mathcal{F}, \mathcal{I} ) := 
\sum_{G \in  \mathcal{F} } s^{\mbox{dord}(G)} +   \sum_{G \in \mathcal{I}} \frac{s^{\mbox{dord}(G)}}{ ( \mbox{free}(G) - 1)! }   \left(  \frac{d}{d s}  \right)^{\mbox{free}(G) - 1} (1 - s)^{-1}
\end{equation}
In the above formula $\mbox{free}(G)$ is the number of free variables in the right hand side of the infinite data set $\mathcal{I}$. \index{$\mbox{free}$, The number of free variables in $\mathcal{I}$ set} So for our example with $G =  u (x_0, y) = f(y)$ we get $\mbox{free}(  u (x_0, y) = f(y) )=1$ and $\mbox{dord}(G) = 0$.  Then the formula yields as before
\begin{equation}
\label{HFcalc}
\mbox{HF} (\mathcal{F}, \mathcal{I} ) :=  s + \frac{1}{1 - s}
\end{equation}
For leading linear {\DPS}, in {\rif}-form with respect to an orderly Riquier ranking, the Differential Hilbert Series gives coordinate independent dimension information.  If a ranking is not orderly then the Differential Hilbert Series is no longer invariant.  For example just consider the initial data for $v_{xx} = v_t$ in an orderly ranking compared to $v_t = v_{xx}$
in a non-orderly ranking.

The output of the \index{{\ID}, Initial data} initial data which partitions the parametric derivatives into disjoint cones of various dimensions $\leq n$, can express this in the rational function form $\mbox{HS}(s) = \frac{P(s)}{(1 - s)^{\textbf{d}}}$
where 
$\textbf{d} = \textbf{d}(\Sys)$ is the differential dimension of $\Sys$. \index{{\textbf{d}}, Differential dimension}  It corresponds to the maximum number of free independent 
variables appearing in the functions for the initial data.
For further information on differential Hilbert Series see \cite{MLPK13:Hilbert}.
The algorithms are simple modifications of those for Gr\"obner bases for modules.

\section{Symmetries \& Mapping Equations}
\label{sec:PreMapEqs}

\subsection{Symmetries}

Infinitesimal Lie {\em point} symmetries for $\Sys$ are found by seeking
vector fields
\begin{equation}
	\label{eq:symmOp}
	V = \sum_{i=1}^n \xi^i(x,u) \dd{x^i} + \sum_{j=1}^m \eta^j(x,u)\dd{u^j}
\end{equation}
whose associated one-parameter group of transformations
\begin{align}
	\label{eq:LieInfTrans}
	{x^*}   &= x +  \xi(x,u) \epsilon + O(\epsilon^2) \nonumber\\
	{u^*}   &= u + \eta(x,u) \epsilon + O(\epsilon^2)
\end{align}
away from exceptional points preserve the jet locus of such systems,  mapping solutions to solutions.
See \cite{BluKu111:Sym,BluKu112:Sym} for applications.
The {\em infinitesimals} $(\xi^i,\eta^j)$ of a symmetry vector field~\eqref{eq:symmOp} for a system of \DEs are found by solving an associated system of linear homogeneous defining equations $S$ (or determining equations) for the infinitesimals.
The defining system $S$ is derived by a prolongation formula for which numerous computer
implementations exist \cite{Car00:Sym,Che07:GeM,Roc11:SAD}.
Lie's classical theory of groups and their algebras requires local analyticity in its defining equations.  Such local analyticity will be a key assumption throughout our paper.

The resulting vector space of vector fields is closed under its commutator.
The commutator of two vector fields for vector fields 
$X  =  \sum_{i=1}^{m+n} \nu^i \dd{z^i}$, $Y =  \sum_{i=1}^{m+n} \mu^i \dd{z^i} $  in a Lie algebra $\mathcal{L}
$ and $z = (x,u)$, is:
\begin{equation}
	\label{eq:com}
	\left[X ,  Y \right] =  X Y - Y X =  \sum_{i=1}^{m+n} \omega^i \dd{z^k}
\end{equation}
where $\omega^k = \sum_{i=1}^{m+n} \left( \nu^i \mu^k_{z^i} - \mu^i \nu^k_{z^i} \right)$.

Similarly, we suppose that the \Target admits symmetry vector fields

\begin{equation}
	\label{eq:vf}
	\hat{V} = \sum_{i=1}^n \hat{\xi}^i(\hat{x},\hat{u}) \dd{\hat{x}^i} + \sum_{j=1}^m \hat{\eta}^j(\hat{x},\hat{u})\dd{\hat{u}^j}
\end{equation}
in the \Target infinitesimals $(\hat{\xi}, \hat{\eta})$ that satisfies a linear homogeneous 
defining system $\hat{S}$ generating a Lie algebra $\hat{\mathcal{L}}$.
Computations with defining systems will be essential in our approach and are
implemented using Huang and Lisle's powerful object \index{{\LAVF}, LieAlgebrasOfVectorFields Maple package} oriented {\LAVF}, {\Maple} package~\cite{LisH:Alg}.

\begin{example}
	\label{ex:ODE3(a)}
	Consider as a simple example the third order nonlinear {\ODE} which is in \rif-form with respect to an orderly ranking
	\begin{equation}
		\label{ODE3}
		u_{{xxx}}  =   {\frac { 3 \left( uu_{{xx}}+{u_{{x}}}^{2}+1 \right) ^{2}}{u(uu_{{x}}+x)}}
		-  {\frac{3 u_{{x}}u_{{xx}}}{u}}
		+ {\frac {8 x \left( uu_{{x}}+x \right) ^{4} \left( {u}^{2}+{x}^{2}+1 \right) }{u({u}^{2}+{x}^{2})}}
	\end{equation}
	at points $u \not = 0, uu_{{x}}+x \not = 0,  {u}^{2}+{x}^{2} \not = 0$.
	When {\Maple}'s \maple{dsolve} is applied to (\ref{ODE3}) it yields no result.  Later in this section, we will discover
	important information about (\ref{ODE3}) using symmetry aided mappings.  It will be used as a simple running example
	to illustrate the techniques of the article. 
	
	The defining system for Lie point symmetries of form
	$\xi(x,u) \dd{x} + \eta(x,u) \dd{u}$  of  (\ref{ODE3}) has {\rif} form with respect to an orderly ranking 
	given by:
	\begin{align}
		\label{eq:ODE3DetSys}
		S  = [& \xi  = -{\frac {\eta\,u}{x}}, \quad
		\eta_{{x,u}}={\frac { \left( u-x \right) 
				\left( u+x \right) \eta}{{u}^{3}x}}+{\frac { \left( {u}^{2}+{x}^{2}
				\right) \eta_{{u}}}{x{u}^{2}}}-{\frac {\eta_{{x}}}{u}}+{\frac {x\eta_
				{{u,u}}}{u}},   \nonumber \\ \qquad
		&	\eta_{{x,x}}=-{\frac { \left( 2\,{u}^{4}-{x}^{2}{u}^{2}+{
					x}^{4} \right) \eta}{{u}^{4}{x}^{2}}}+{\frac { \left( {u}^{2}+{x}^{2}
				\right) \eta_{{u}}}{{u}^{3}}}+2\,{\frac {\eta_{{x}}}{x}}+{\frac {{x}^
				{2}\eta_{{u,u}}}{{u}^{2}}},    \\
		&	\eta_{{u,u,u}}=-{\frac { \left( 16\,{u}^{8}
				+24\,{u}^{6}{x}^{2}+8\,{u}^{4}{x}^{4}+16\,{u}^{6}+8\,{u}^{4}{x}^{2}+3
				\,{u}^{2}+3\,{x}^{2} \right) \eta}{ \left( {u}^{2}+{x}^{2} \right) {u}
				^{3}}}    \nonumber \\
		& -{\frac { \left( 8\,{u}^{6}{x}^{2}+8\,{u}^{4}{x}^{4}+8\,{u}^{4}{
					x}^{2}-3\,{u}^{2}-3\,{x}^{2} \right) \eta_{{u}}}{{u}^{2} \left( {u}^{2
				}+{x}^{2} \right) }} +8\,{\frac {{u}^{3}x \left( {u}^{2}+{x}^{2}+1
				\right) \eta_{{x}}}{{u}^{2}+{x}^{2}}} \nonumber ]
	\end{align}
	{\todo{[R3C1]}}Its corresponding initial data is
	\begin{align}
		\hspace{-0.5cm} \mbox{\ID} (S) = [\eta \left( x_{{0}},u_{{0}} \right) =c_{{1}},  \eta_{{x}} \left( x_{{0}
		},u_{{0}} \right) =c_{{2}},  \eta_{{u}} \left( x_{{0}},u_{{0}} \right) =
		c_{{3}},  \eta_{{u,u}}   \left( x_{{0}},u_{{0}} \right) =
		c_{{4}}]
	\end{align}
	There are $4$ arbitrary constants in the initial data at regular points $(x_0, u_0)$,
	so (\ref{ODE3}) has a $4$ dimensional local Lie algebra of symmetries $\LieAlg$ in a neighborhood of such points:
	$\dim \mathcal{L} = 4$.
	The structure of $\LieAlg$ of (\ref{eq:ODE3DetSys}) can be algorithmically determined without integrating the defining system \cite{RLB92:Alg,LisleReidInfinite98,LisH:Alg,Huang:Thesis}:
	\begin{align}
		[Y_{{1}},Y_{{2}}] &=-Y_{{1}}-2\,Y_{{2}}, \quad [Y_{{1}},Y_{{3
		}}]=Y_{{1}}-2\,Y_{{3}},\quad [Y_{{1}}, Y_{{4}}]=-2\,Y_{{4}}, \nonumber  \\
		[Y_{{2}},Y_{{3}}] &=Y_{{2}}+Y_{{3}}, \quad [Y_{{2}},Y_{{4}}]=Y_{{4}},\quad [Y_{{3}},Y_{{4}}]=-Y_{{4}}
	\end{align}
	where a regular point ($x_0 = 1, u_0 = 1$) was substituted into the relations.
\end{example}
But what can such symmetry information tell us about nonlinear systems $\Sys$ such as the above {\ODE} using mappings?
In particular, in this paper we focus on the question of when a system $\Sys$ can be mapped to a linear system $\Syshat$.
Throughout this paper we maintain blanket local analyticity assumptions.
So the case of a single differential equation $\Syshat$ has the form $\mathcal{H} \hat{u} = f(\hat{x})$ where $\mathcal{H}$ is a linear differential operator with coefficients that are analytic functions of $\hat{x}$ and $f$ is also 
analytic.  Lewy's famous counterexample of a single linear differential equation in 3 variables, of order 1, where 
$\mathcal{H}$ is analytic with smooth inhomogeneous term, without smooth solutions, provides a counterexample in the smooth case. 
Then supposing we have the existence of a local analytic solution $\tilde{u}$ in a neighborhood of $\hat{x}_0$, 
in this neighborhood the point transformation $\hat{u} \rightarrow  \hat{u} - \tilde{u}$ implies that without loss 
we can consider 
$\Syshat$ to be a homogeneous linear differential equation $\mathcal{H} \hat{u} = 0$ where 
$\hat{u}  \in \mathcal{A}(\hat{x}_0, \delta)$, the set of analytic functions on some sufficiently small disk
$| \hat{x} - \hat{x}_0 | < \delta$.
Then solutions of $\Syshat$ satisfy the superposition property
$\mathcal{H} (\hat{v}  + \hat{w}) =  \mathcal{H} \hat{v}  +  \mathcal{H}  \hat{w}  = 0$.
This corresponds to point symmetries generated by the Lie algebra of vectorfields
\begin{equation}
	\label{L*}
	\hat{\mathcal{L}}^* := \left\{ \hat{v}(\hat{x}) \frac{\partial}{\partial \hat{u}}: \mathcal{H} \hat{v}(\hat{x}) = 0 \; and \; 
	\hat{v}  \in \mathcal{A}(\hat{x}_0, \delta)  \right\}
\end{equation}
Consequently, assuming the existence of a local analytic map, $\dim \hat{\mathcal{L}}^* = \dim \Syshat = \dim \Sys$.
If $\Sys$ is an {\ODE} of order $d \geq 2$ then $\dim \hat{\mathcal{L}}^* = \dim \Syshat = \dim \Sys = d$.
Similarly the superposition property
$\mathcal{H} (c\hat{v}) = c \mathcal{H} \hat{v}  = 0$
corresponds to a $1$ parameter family of scalings with symmetry vectorfield $\hat{u} \frac{\partial}{\partial \hat{u}}$.
So we get the well-known and obvious result that an {\ODE} of order $d$ that can be mapped 
to a linear {\ODE} must have 
$\dim {\mathcal{L}} = \dim \hat{\mathcal{L}} \geq d + 1$.
Similarly if $\Sys$ is linearizable and $\dim \Sys = \infty$ then $\dim {\mathcal{L}} = \dim \hat{\mathcal{L}} = \infty$ with similar properties for systems. 
The Lie sub-algebra $\hat{\mathcal{L}}^*$ is easily shown to be `abelian` by direct computation of commutator \cite{PetOlv107:Sym} in both the finite and infinite dimensional case.
Indeed consider the so-called derived algebra $\mathcal{L}' = \mbox{DerivedAlgebra}(\mathcal{L})$, which is the Lie subalgebra of 
$\mathcal{L}$ generated by commutators of members of $\mathcal{L}$ and similarly for $\hat{\mathcal{L}}'$.
By direct computation of commutators $\hat{\mathcal{L}}^*$ is a sub-algebra of $\hat{\mathcal{L}}'$ (e.g. $[ \hat{v}(\hat{x}) \frac{\partial}{\partial \hat{u}},  
\hat{u} \frac{\partial}{\partial \hat{u}} ] = \hat{v}(\hat{x}) \frac{\partial}{\partial \hat{u}} \in \hat{\mathcal{L}}$).
Thus, a necessary condition for the existence of a map $\Psi$ to a linear target is that 
$\hat{\mathcal{L}}'$ has a $d$ dimensional abelian subalgebra in the finite and infinite dimensional
cases (see Olver \cite{PetOlv107:Sym}).
{\todo{[R1C7]}}
In the preceding paragraph we considered the case of a single dependent variable, which is easily extended to the multivariate
case.  For example in equation  (\ref{L*}) the symmetry generator $\hat{v}(\hat{x}) \frac{\partial}{\partial \hat{u}}$ can be
replaced by $ \sum_{i = 1}^{m} \hat{v}^i (\hat{x}) \frac{\partial}{\partial \hat{u}^i} $ for the case of a system.

{\todo{[R4C4]}}
Lyakhov, Gerdt and Michels \cite{LGM101:LG} use this to implement a remarkable algorithm to determine the existence of a linearization
for a single {\ODE} of order $d$.  See Algorithm \ref{alg:LGMLinTest}
and \cite{LGM101:LG,ML102:LSA} for further background.
There are two main cases.  The first is when the nonlinear {\ODE} has a Lie symmetry algebra of maximal dimension, as shown in 
Step 5 of Algorithm \ref{alg:LGMLinTest}.  Such maximal cases are always linearizable. These occur for $d = 1, 2$ where the maximal dimensions of $\mathcal{L}$ are 
$\infty$ and $8$ respectively, and for $d > 2$ where the maximal dimension of $\mathcal{L}$ is $d+4$.
The second main case is sub-maximal and occurs for $d > 2$ when $\dim \mathcal{L} = d+1$ or $\dim \mathcal{L} = d+2$. 
\index{$\mbox{DetSys}$, Determining system}

\begin{algorithm}[h!]
	\caption{ ${\mbox{LGMLinTest}}(\Sys)$   } \index{${\mbox{LGMLinTest}}$, Algorithm\ref{alg:LGMLinTest}} 
	\begin{algorithmic}[1]
		\Statex {\bf Input}: a leading linear {\ODE} $\Sys$ solved for its highest derivative of order $\geq 1$
		\Statex {\bf Output}:  Lin = true if $R$ linearizable otherwise Lin = false
		\State	Lin:= false  		                                                          
		\State Compute $S := \mbox{ThomasDecomposition}(\mbox{DetSys}(R))$ 
		\State Find $\dim\mathcal{L} := \dim S$, $d := \mbox{difforder}(R)$                       
		\State ComRels := Structure$(S)$   {\todo{[R1C8]}}                                            
		\State	\textbf{if} $d=1$ or ($d=2$ and $\dim \mathcal{L}=8$) or ($d>2$ and $\dim \mathcal{L}= d+4$)   
		\textbf{then} Lin := true         
		\Statex	 \textbf{else if} $d>2$ and ($\dim \mathcal{L}= d+1$ or $\dim \mathcal{L}= d+2$) \textbf{then} 
		\Statex	\hskip12pt   $\mathcal{L}'  := \mbox{DerivedAlgebra}(\mbox{ComRels})$
		\Statex	\hskip12pt \textbf{if} IsAbelian$(\mathcal{L}')$ and $d = \dim(\mathcal{L}')$ \textbf{then} Lin := true  
		\Statex \hskip12pt \textbf{end if} 
		\Statex  \textbf{end if}
		\State	\textbf{return} Lin
	\end{algorithmic}
	\label{alg:LGMLinTest}
\end{algorithm}

\begin{example}
	\label{ex:ODE3(b)}
	We illustrate the above discussion and Algorithm \ref{alg:LGMLinTest} by a continuation of Example \ref{ex:ODE3(a)}. For that example $d = 3$ and $\dim \mathcal{L} = 4$.
	Then from the commutation relations the derived algebra $\mathcal{L}'$ is generated by 
	\begin{equation}
		[Z_{{1}}=Y_{{1}}-2\,Y_{{3}},\; Z_{{2}}=Y_{{2}}+Y_{{3}},\;Z_{{3}}=Y_{{4}}]
	\end{equation}
	Thus $\dim \mbox{DerivedAlgebra}( \mathcal{L}) = \dim \mathcal{L}' = 3$. Also its structure is easily found as 
	\begin{equation}
		[Z_{{1}},Z_{{2}}] = [Z_{{1}},Z_{{3}}] = [Z_{{2}},Z_{{3}}] = 0
	\end{equation}
	{\todo{[R4C6]}}
	so $ \mathcal{L}'$ is abelian and by Algorithm  \ref{alg:LGMLinTest}, (\ref{ODE3}) is exactly linearizable.

\end{example}

{\todo{[R4C6]}} The Algorithms introduced by Lyakhov et al. \cite{LGM101:LG} have two stages: the first given above is to determine whether the system is lineaizable.  The second 
stage is to attempt to construct an explicit form for the mapping by integration. 
A fundamental algorithmic tool for both stages is the \maple{ThomasDecomposition} algorithm
which is a differential elimination algorithm {\todo{[R1C9], [R4C5]}} which outputs a disjoint decomposition of 
a {\DPS} finer than that of 
\cite{Bou95:Rep} or \cite{Rus99:Exi}.  It is based on the work of Thomas
\cite{Thomas106:Tdec,Robertz106:Tdec}.
The algorithm
is available in distributed Maple 18 and later versions. We note that the construction step involves heuristic integration.
The algorithm that they use to construct a system for the mapping to a linear {\ODE}, before it is reduced using \maple{ThomasDecomposition}, is related to
the algorithm \maple{EquivDetSys} given in our introductory paper \cite{MohReiHua19:Intro}.
It is expensive as we illustrate later with examples. One of the contributions of our paper is to find a potentially more efficient algorithm that avoids the application of the full nonlinear equivalence equations generated by \maple{EquivDetSys} in \cite{MohReiHua19:Intro}.

\subsection{Bluman-Kumei Mapping Equations}

Assume the existence of a local analytic invertible map $\Psi = (\psi, \phi)$ between the \Source system $\Sys$ and the \Target system $\Syshat$, with Lie symmetry algebras $\mathcal{L}$, $\hat{\mathcal{L}}$ respectively.  Applying $\Psi$ to the infinitesimals $(\hat{\xi}, \hat{\eta})$ of a vectorfield in
$\hat{\mathcal{L}}$ yields what we will call \index{$M_{\mbox{\sc{BK}}}$, Bluman-Kumei mapping equations} the \textit{Bluman-Kumei (BK) mapping equations}:
\begin{equation}
	\label{eq:BK}
	{\small{M_{\mbox{\sc{BK}}}(\mathcal{L}, \hat{\mathcal{L}}) }} =
	\left\{ \begin{array}{ccl}
		\hat{\xi}^k (\hat{x}, \hat{u})  & = & \sum_{i=1}^n \xi^i(x,u)  \frac{\partial \psi^k}{\partial x^i} + \sum_{j=1}^m \eta^j(x,u)\frac{\partial \psi^k}{\partial u^j}    \vspace{0.5cm}  \\
		\hat{\eta}^\ell (\hat{x}, \hat{u})  & = & \sum_{i=1}^n \xi^i(x,u) \frac{\partial \phi^\ell}{\partial x^i}+\sum_{j=1}^m \eta^j(x,u)\frac{\partial \phi^\ell}{\partial u^j}   \\
	\end{array}\right. 
\end{equation}
where $1 \leq k \leq n$ and $1 \leq \ell \leq m$, and $(\xi, \eta)$ are infinitesimals of Lie symmetry vectorfields in $\mathcal{L}$.
See Bluman and Kumei \cite{BK105:BluKu,Blu10:App} for details and generalizations (e.g. to contact transformations). Note that all quantities on the RHS of the BK mapping equations~\eqref{eq:BK} are functions of $(x,u)$ including $\phi$ and $\psi$. See \cite[Example 1]{MohReiHua19:Intro} for an introductory example of mappings and the
examples in \cite{Blu10:App}. %

\begin{remark} When considered together with $\hat{x} = \psi(x,u), \hat{u} = \phi(x,u)$ the BK mapping equations ~\eqref{eq:BK} are a change of variables from $(x,u)$ to $(\hat{x}, \hat{u})$ coordinates.  Simply interchanging target and source variables then yields the inverse of the BK mapping equations below.  Considered together with $x = \hat{\psi}(\hat{x}, \hat{u}), u = \hat{\phi}(\hat{x}, \hat{u})$ these are 
a change of variables from $(\hat{x}, \hat{u})$ to $(x,u)$ coordinates. 
\label{rem:InvBK}
\begin{equation}
	\label{eq:InvBK}
	{\small{M_{\mbox{\sc{BK}}}( \hat{\mathcal{L}}, \mathcal{L}) }} =
	\left\{ \begin{array}{ccl}
		\xi^k(x,u)  & = & \sum_{i=1}^n \hat{\xi}^i (\hat{x}, \hat{u})  \frac{\partial \hat{\psi}^k}{\partial \hat{x}^i} + \sum_{j=1}^m \hat{\eta}^j (\hat{x}, \hat{u}) \frac{\partial \hat{\psi}^k}{\partial \hat{u}^j}    \vspace{0.5cm}  \\
	    \eta^\ell (x,u)  & = & \sum_{i=1}^n \hat{\xi}^i (\hat{x}, \hat{u}) \frac{\partial \hat{\phi}^\ell}{\partial \hat{x}^i}+\sum_{j=1}^m\hat{\eta}^j (\hat{x}, \hat{u}) \frac{\partial \hat{\phi}^\ell}{\partial \hat{u}^j}   \\
	\end{array}\right. 
\end{equation}
We note the following relation between Jacobians in $(x,u)$ and $(\hat{x},\hat{u})$ coordinate systems
\begin{equation}
	\label{eq:PsiPsihat}
	\frac{\partial(\hat{\psi},\hat{\phi})}{\partial(\hat{x}, \hat{u})} =  \left[ \frac{\partial(\psi,\phi)}{\partial(x, u)} \right]^{-1} 
\end{equation}

\end{remark}

If an invertible map $\Psi$ exists mapping $\Sys$ to $\Syshat$ then it most generally 
depends on $\dim(\mathcal{L}) = \dim(\hat{\mathcal{L}})$ parameters.
But we only need one such $\Psi$.  So reducing the number of such parameters, e.g.,\ by restricting to 
a Lie subalgebra $\mathcal{L}^{'}$ of $\mathcal{L}$ with corresponding Lie subalgebra $\hat{\mathcal{L}}^{'}$ of $\hat{\mathcal{L}}$
that still enables the existence of such a $\Psi$, is important in reducing the computational
difficulty of such methods.  
We will use the notation $S',\hat{S}'$ to denote the symmetry defining systems of  Lie sub-algebras $\mathcal{L}^{'} $,
$\hat{\mathcal{L}}^{'}$ respectively.
See \cite{Blu10:App,PetOlv107:Sym} 
for discussion on this matter.

For mapping from nonlinear to linear systems, a natural candidate for $\mathcal{L}^{'}$ is the
$ \mbox{DerivedAlgebra}  (\mathcal{L}) $, and the natural target Lie symmetry algebra is $\hat{\mathcal{L}}^{*}$ corresponding to the 
superposition defined in (\ref{L*}).

\begin{example}
	\label{ex:ODE3(c)}
	This is a continuation of Examples \ref{ex:ODE3(a)} and \ref{ex:ODE3(b)} concerning (\ref{ODE3}).  Here we will use the BK mapping equations~\eqref{eq:BK}
	where 
	$\mathcal{L}^{'} = \mbox{DerivedAlgebra}  (\mathcal{L}) $.
	
	For the construction of $\Psi$ we actually need differential equations for $ \mathcal{L}'$, in
	addition its structure which are provided by the algorithm $\maple{DerivedAlgebra}$
	in the {\LAVF} package which for (\ref{ODE3}) yields its {\rif}-form:
	\begin{align}
		\label{eq:ODE3_L'}
		S'  = [\xi = &-{\frac {\eta\,u}{x}},\eta_{{x}}={\frac { \left( {u}^{2}+{x}^{2}
				\right) \eta}{x{u}^{2}}}+{\frac {\eta_{{u}}x}{u}},    \\
		\eta_{{u,u,u}}&=-{
			\frac { \left( 8\,{u}^{8}+8\,{u}^{6}{x}^{2}+8\,{u}^{6}+3\,{u}^{2}+3\,{
					x}^{2} \right) \eta}{ \left( {u}^{2}+{x}^{2} \right) {u}^{3}}}+3\,{
			\frac {\eta_{{u}}}{{u}^{2}}}]  \nonumber
	\end{align}
	The derived algebra is then shown by {\LAVF} commands to be both $3$ dimensional and abelian.  Moreover, its determining system (\ref{eq:ODE3_L'}) is much simpler than the determining system of $\mathcal{L}$ given in (\ref{eq:ODE3DetSys}).
	Crucially it means we can exploit this determining system using the BK mapping equations.  Since the target infinitesimal generator is $\hat{\xi} \dd{\hat{x}} + \hat{\eta} \dd{\hat{u}} = 0 \cdot \dd{\hat{x}} + \hat{\eta} (\hat{x}) \dd{\hat{u}}$ 
	where $\mathcal{H} \hat{u}  = 0$.  So  $ \hat{\xi} = 0$ and $\mathcal{H} \hat{u}  = 0$ and the BK equations are:
	\begin{equation}
		\label{eq:BKODE3}
		{\small{M_{\mbox{\sc{BK}}}(\mathcal{L}^{'}, \hat{\mathcal{L}}^{*})}} =
		\left\{ \begin{array}{ccl}
			0  		& = &  \xi (x,u)  \frac{\partial \psi}{\partial x} +  \eta(x,u)\frac{\partial \psi}{\partial u}    \vspace{0.5cm}  \\
			\hat{\eta} (\hat{x}, \hat{u})  & = &  \xi(x,u) \frac{\partial \phi}{\partial x}+\eta(x,u)\frac{\partial \phi}{\partial u}   \\
		\end{array}\right. 
	\end{equation}
	{\todo{[R4C7]}}
	where $\hat{\eta}_{\hat{u}} = 0$ and $\mathcal{H} \hat{u}  = 0$.
	These equations are an important necessary condition for the linearization of (\ref{ODE3}) and this will be exploited in Section \S \ref{sec:Algorithms} in the computation of the mapping.

\end{example}

\section{Algorithms and Preliminaries for the MapDE Algorithm}
\label{sec:Algorithms}

The \MapDE algorithm introduced in \cite{MohReiHua19:Intro} is extended to determine if there exists a mapping of a nonlinear source $\Sys$ to some linear target $\Syshat$, using the target input option $\maple{\Target}= \maple{LinearDE}$.

\subsection{Symmetries of the linear target and the derived algebra}

\label{sec:SymDer}

We summarize and generalize some aspects of the discussion in \S \ref{sec:Intro} and \S \ref{sec:PreMapEqs}.
The following theorem is a straightforward consequence of the necessary conditions in  \cite{Blu10:App}
where we have also required that the target system is in {\rif}-form.

\begin{theorem}[Superposition symmetry for linearizable systems]
	\label{thm:SupSymLinSys}
	Suppose that the analytic system $\Sys$ is exactly linearizable by a local holomorphic diffeomorphism
	$\hat{x} = \psi(x,u), \hat{u} = \phi(x,u)$ to yield a linear target system.  Then $\Syshat$ locally takes the form
	\begin{equation}
		\label{eq:Hsy}
		\Syshat: \mathcal{H} \hat{u}(\hat{x}) = 0
	\end{equation}
	where $\mathcal{H}$ is a vector partial differential operator, with coefficients that are local analytic functions of $\hat{x}$ and {\todo{[R4C8]}} the system (\ref{eq:Hsy}) is in {\rif}-form  with respect to an orderly ranking.  Moreover $\Syshat$ admits the symmetry vector field $\sum_{j=1}^m \hat{\eta}^j (\hat{x}) \frac{\partial}{\partial \hat{u}^j}$:
	\begin{equation}
		\label{eq:Hsym}
		\hat{S}^* := \left\{ \hat{\xi}^i = 0, \; \mathcal{H} \hat{\eta} = 0,  \;  \hat{\eta}^j_{\hat{u}^k}  = 0: 1 \leq i \leq n,  \; 1 \leq j,k \leq m \right \}.
	\end{equation}
\end{theorem}

From the previous discussion, computation of determining systems for derived algebras is important in both the finite and 
infinite cases.  In the Remark below we sketch what appears to be the first algorithm to compute such systems in the infinite case.

\begin{remark}[Algorithm for computation of infinite Derived Algebras]
\label{rem:DerivedAlgebraAlgorithm}
A simple consequence of the commutator formula (\ref{eq:com}) is that the commutators generate a Lie 
algebra which is called the derived algebra.
Lisle and Huang \cite{LisH:Alg} implement efficient algorithm in the {\LAVF} package to compute the determining system for the derived algebra for finite dimensional Lie algebras of vectorfields.
We have made a first implementation in the infinite dimensional case, 
together with the Lie pseudogroup structure relations \cite{RLB92:Alg,LisleReidInfinite98}.   First each of the $\nu$, $\mu$, $\omega$ in the commutation relations (\ref{eq:com}) must satisfy the determining system
of $\mathcal{L}$ so we enter three copies of those determining systems.  We then reduce the combined system
using a block elimination ranking which ranks any derivative of $\omega$ strictly less than those of $\mu, \nu$.
The resulting block elimination system for $\omega$ generates the derived algebra in the infinite case.
\end{remark}

The commutator between any superposition generator and the scaling symmetry admitted by linear systems yields
\begin{equation} 
	\label{eq:DerCom}
	\left[ \sum_{i = 1}^{m}  \hat{v}^i(\hat{x}) \frac{\partial}{\partial \hat{u}^i} \; , \; \sum_{i = 1}^{m} \hat{u}^i \frac{\partial}{\partial \hat{u}^i}  \right] =  \sum_{i = 1}^{m} \hat{v}^i (\hat{x}) \frac{\partial}{\partial \hat{u}^i}
\end{equation}
So we have the following result as an easy consequence (See Olver \cite{PetOlv107:Sym} for related discussion in both the finite and infinite case). 

\begin{theorem}
	\label{thm:DerivedAlgebraThm}
	Suppose that the analytic system $\Sys$ is exactly linearizable by a local holomorphic diffeomorphism
	$\hat{x} = \psi(x,u), \hat{u} = \phi(x,u)$, to yield a linear target system (\ref{eq:Hsym})
	and $\mathcal{L}$, $\mathcal{L}'$ are the Lie symmetry algebra and its derived algebra for $\Sys$.
	Also, let $\hat{\mathcal{L}}$, $\mathcal{\hat{L}}'$ be the corresponding algebras for ${\Syshat}$.
	Let $\mathcal{L}^*$, $\mathcal{\hat{L}}^*$ be the superposition algebras under $\Psi$.
	Then $\mathcal{\hat{L}}^*$ is a subalgebra of $\mathcal{\hat{L}}'$ and
	$\mathcal{L}^*$ is a subalgebra of $\mathcal{L}'$.
	Moreover $\mathcal{L}^*$ and $\mathcal{\hat{L}}^*$ are abelian.
\end{theorem}
We wish to determine if a system $\Sys$ is linearizable and if so, characterize the target $\Syshat$, i.e 
$\mathcal{H} \hat{\eta}= 0$.   But initially we don't know $\mathcal{H}$. 
One approach is to write a general form for this system that specifies $\hat{S}^*$ 
with undetermined coefficient functions 
whose form is established in further computation.  See for example, \cite{LGM101:LG} use this approach in the case of a single {\ODE}, but don't consider the Bluman-Kumei mapping system.{\todo{[R4C9]}} We will apply our method to a test set of {\ODE} (\ref{LGMTestODE}) given in \cite{LGM101:LG}.
Instead, we only include $\xi^i = 0, \hat{\eta}^j_{\hat{u}^k}  = 0$ and don't include $\mathcal{H} \hat{\eta} = 0$.
Thus, we only include a subset $\hat{S}^\star$ of $\hat{S}^*$, denoting the truncated system as 
{\todo{[R4C10]}}
\begin{equation}
	\label{eq:Sstar}
	\hat{S}^\star :=  \left\{ \hat{\xi}^i = 0, \hat{\eta}^j_{\hat{u}^k}  = 0:  1 \leq i \leq n,  \; 1 \leq j,k \leq m  \right\} 
\end{equation}
and allow $\mathcal{H} \hat{\eta} = 0$ to be found naturally later in the algorithm.
We note that $\hat{S}^\star$ are the defining equations of a (usually infinite) Lie pseudogroup.

\subsection{Algorithm PreEquivTest  for excluding obvious nonlinearizable cases}

\label{sec:PreEquivTest}

As discussed in \S \ref{sec:PreMapEqs}, $\Sys$ being linearizable implies that the superposition is in its Lie symmetry algebra and is a coordinate change of $\Sys$.  This implies {\todo{[R4C11]}} some fairly well-known efficient tests for screening out obvious non-linearizable cases.  For linearization  necessarily $\dim{S} \geq d+1$ and $\dim{S'} \geq d$ for finite $d$. For $d = \infty$, necessarily $\dim{S} = \infty = \dim{S'}$ and in terms of differential dimensions $\textbf{d}(S) \geq \textbf{d}(S') \geq \textbf{d}(R)$.  Note that Algorithm \ref{Alg:PreEquivTest} returns null, if all its tests are true.
The most well-known of the above tests occur when $d = \infty$ and $\dim{S} = \infty$ and are given for example in 
Bluman and Kumei \cite{BluKu112:Sym}.  Also see Theorem 6.46 in Chapter 6 of Olver \cite{PetOlv107:Sym}.

\begin{algorithm}[h!]
	\caption{\PreEquivTest }  \index{{\PreEquivTest}, Algorithm \ref{Alg:PreEquivTest}} 
	\begin{algorithmic}[1] 
		\Statex ${\mbox{\PreEquivTest}}(\Sys, \mbox{IDR}, \mbox{\sc{IDS}}, \mbox{IDS}')$ 
		\Statex {\bf Input}:  $\Sys$ is a leading linear {\DPS} system in {\dec}-form $\Sys$
							with no leaders of order $0$ 
		\Statex \hskip24pt   with respect to an orderly Riquier ranking
		\Statex \hskip12pt $\mbox{IDR}$, $\mbox{IDS}$, $\mbox{IDS}'$ are respectively the initial data for $\Sys$, $S$ and $S'$
		\Statex \hskip24pt   where $S$ and $S'$ are the symmetry determining systems for $\mathcal{L}$ and $\mathcal{L}'$
		\Statex {\bf Output}:  [IsLinearizable, DimInfo] 
		\Statex \hskip24pt  IsLinearizable = false if one of the necessary conditions $T_j$ tests false
		\Statex \hskip24pt  DimInfo is the dimension info (dimension, differential dimension and
		\Statex \hskip24pt  Differential Hilbert Function, computed for each of $R, S, S'$.
		\State Set IsLinearizable := null. Apply DifferentialHilbertFunction to IDR, IDS, IDS' to get 
		\State  \hskip12pt	DimInfo := [ $\dim{\Sys}$, $\textbf{d}(R)$, HF(R), $ \dim{S}$, $\textbf{d}(S)$, HF(S) $, \dim{S'}$, $\textbf{d}(S')$, HF$(S')$]
		\State		\textbf{if} $d = \infty$ \textbf{then} 
		\Statex		\hskip12pt 			$T_1$ := evalb($\dim{S} < \infty$)
		\Statex		\hskip12pt 			$T_2$ := evalb($\dim{S'}< \infty$)
		\Statex		\hskip12pt 			$T_3$ := evalb($\textbf{d}(S)  < \textbf{d}(R)$)
		\Statex		\hskip12pt 			$T_4$ := evalb($\textbf{d}(S') < \textbf{d}(R)$)
		\Statex		\textbf{else if} $d < \infty$ \textbf{then}
		\Statex		\hskip12pt 			$T_5$ := evalb($\dim{S} < d+1$)
		\Statex		\hskip12pt 			$T_6$ := evalb($\dim{S'} < d$)
		\Statex		\textbf{end if}  					
		\State 		\textbf{if}  $\land_{i = 1}^{i = 6} T_i $ = false \textbf{then} IsLinearizable := false;  \textbf{end if}  
		\Statex		\textbf{return} [ IsLinearizable, DimInfo ]
		\Statex		  
	\end{algorithmic}
	\label{Alg:PreEquivTest}
\end{algorithm}

\subsection{Algorithm ExtractTarget for extracting the linear target system} 
\index{\text{ExtractTarget}, Algorithm \ref{Alg:ExtractTarget}} 
\label{Alg:ExtractTarget}

When a system $\Sys$ is determined to be linearizable by Algorithm \ref{Alg:MapDE2Linear}, the conditions
for linearizability will yield a list of cases $\bigcup_{c \in C} Q_c$ where each $Q_c$ is in {\rif}-form.
 To implicitly determine the target linear system $\Syshat$ for a case $Q_c$, Algorithm $\text{ExtractTarget}$
is applied to $Q_c$.  It first selects from $Q_c$ the linear homogeneous differential sub-system $ R_c^*$ in $\xi(x,u)$, $\eta(x,u)$ with coefficients depending on $(x,u,\psi, \phi)$ in the $(x,u)$ coordinates.  
Algorithm $\text{ExtractTarget}$ then applies the inverse BK transformations (\ref{eq:InvBK}) to convert the system $ R_c^*$ to $(\hat{x},\hat{u})$ coordinates, after which $\hat{\xi} = 0$ and 
also $\hat{\eta}(\hat{x},\hat{u}) =\hat{ \eta}( \hat{x} )$ is imposed.
This yields $ R_c^*$ as a system $\hat{ R}_c^*$ which is a linear homogeneous differential system in $\hat{\eta}(\hat{x})$ with coefficients in $\hat{x}, \hat{u},\hat{\psi}, \hat{\phi}$.
Though $ R_c^*$ is in {\rif}-form, $\hat{ R}_c^*$ is not usually in {\rif}-form so another application of {\rif} is applied, to yield 
$\hat{ R}_c^*$ in {\rif}-form for $\hat{\eta}(\hat{x})$; case splitting is not required here.
As shown in the proof of Algorithm \ref{Alg:MapDE2Linear}, the coefficients of $\hat{ R}_c^*$ depend only on $\hat{x}, \hat{\psi}, \hat{\phi}$ and not on $\hat{u}$. 
The proof of Algorithm \ref{Alg:MapDE2Linear} also shows that $\hat{\eta}(\hat{x})$ can be replaced in $\hat{ R}_c^*$ with $\hat{u}(\hat{x})$ yielding the target linear homogeneous differential equation $\Syshat_c$.

\subsection{Heuristic integration for the Mapping functions in MapDE using PDSolve }

\label{Method:PDSolve}

This routine is still in the early stages of its development.
The heuristic integration routine \maple{PDSolve}, basically a simple interface to {\Maple}'s 
\maple{pdsolve} which is applied to the $\Psi$ sub-system (the sub-system with highest derivatives in $\psi, \phi$) together with its inequations in $Q_c$ to attempt to find an explicit form of the mapping $\Psi$.
We naturally use a block elimination ranking in $\Psi$ sub-system where all derivatives of $\phi$ are higher
than all derivatives of $\psi$. Then we attempt to solve uncoupled subsystem for $\psi$ using
the \LAVF routine Invariants, which depends on integration, and subsequently solve the substituted system 
for $\phi$ by \maple{pdsolve}.

The geometric idea is that in the $(\hat{x}, \hat{u})$ coordinates the Lie symmetry generator corresponding to 
linear superposition has the form $\sum_k \hat{\eta}^k(\hat{x}) \frac{\partial}{\partial \hat{u}^k }$ and generates 
an abelian Lie algebra with obvious invariants $\hat{x}$.  
So the independent variables for the target linear equation are invariants of this vector field acting on the base space of variables $(\hat{x}, \hat{u})$, and thus on $(x,u)$ space via the map $\Psi$.
The process of integrating the mapping equations first starts with the determination of these invariants in terms of 
$(x,u)$ using the {\LAVF} command \maple{Invariants}.  If the integration is successful this yields $\hat{x}^j = \psi^j = I^j(x,u)$, for $j = 1 , \cdots , n$.  Then substitution into the $\Psi$  sub-system, yields a system with dependence 
only on the $\hat{\phi}$ mapping functions, which we attempt to integrate using Maple's \maple{pdsolve}.

	


\section{The MapDE Algorithm}

\label{sec:MapDE2Linear}

The main subject here is  Algorithm \ref{Alg:MapDE2Linear} which makes heavy use of differential-elimination completion ({\dec})  algorithms, which in our current implementation is 
the {\rif} algorithm accessed via {\Maple}'s $\maple{rifsimp}$.
Other {\dec} algorithms could be used such as \maple{ThomasDecomposition}, \maple{RosenfeldGroebner} or \maple{casesplit}.   {\todo{[R1C9]}}

In \S \ref{sec:MapDEcode} we will describe pseudo-code for \MapDE. 
In \S \ref{sec:NotesMapDE2Linear} we will give notes about the steps of \MapDE
and in \S \ref{sec:ProofMapDE} we will a proof of correctness of \MapDE.

\subsection{Pseudo-code for the MapDE algorithm}

\label{sec:MapDEcode}

Here we describe the pseudo-code for \MapDE.


\begin{algorithm}[h!]
	\caption{ \MapDE with \maple{Target}=\maple{LinearDE}  }   \index{{\MapDE}, Algorithm \ref{Alg:MapDE2Linear}} 
	\begin{algorithmic}[1] 
		\Statex ${\mbox{\MapDE}}(\Source, \Target, \Map, Options)$ 
		\Statex {\bf Input}:  
		\Statex \hskip12pt   \Source:  \hskip4pt  A leading linear {\DPS} system in {\dec}-form $\Sys$
		with no leaders of order $0$
		\Statex \hskip55pt   with respect to an orderly Riquier ranking; vars $[x,u]$, $[\xi, \eta]$
		\Statex \hskip12pt   \Target: \hskip6pt \maple{Target}=\maple{LinearDE}
		\Statex \hskip12pt   \Map:    \hskip14pt   $\Psi$
		\Statex \hskip12pt   $Options$:    \hskip2pt   Additional options (for strategies, outputs, etc).  {\todo{[R1C10]}}
		\Statex {\bf Output}: 
		\Statex \hskip12pt   IsLinearizable = false if there $\nexists$ an invertible local linearization $\Psi$
		\Statex \hskip12pt   IsLinearizable = true if there $\exists$ an invertible local linearization $\Psi$ and
		\Statex \hskip20pt   	--- Collection of cases in {\rif}-form yielding such linearizations
		\Statex \hskip20pt   	--- Implicit form of the target linear equation for $\Syshat$ (see \ref{Alg:ExtractTarget} )
		\Statex \hskip20pt   	--- Explicit form of $\Psi$ if the heuristic method \maple{PDSolve} is successful (see \S \ref{Method:PDSolve})
		\State  Set IsLinearizable := null. Compute $ \mbox{IDR} := \mbox{ID}(\Sys)$ 
		\State		Let $\mathcal{L}' = \maple{DerivedAlgebra}(\mathcal{L})$ and compute:							 
		\Statex 		\hskip12pt $S:= \mbox{\dec}(\mbox{DetSys}(\mathcal{L}))$,  $S' := \mbox{\dec}(\mbox{DetSys}(\mathcal{L}'))$ 
		\Statex		\hskip12pt $ \mbox{IDS} :=\mbox{ID}(S )$, $ \mbox{IDS}' :=\mbox{ID}(S')$
		
		\State		[ IsLinearizable, DimInfo ] := $ {\mbox{\PreEquivTest}}(\Sys, \mbox{IDR}, \mbox{\sc{IDS}}, \mbox{IDS}')$
		\Statex		\hskip12pt \textbf{if} IsLinearizable = false \textbf{then} \textbf{return} [ IsLinearizable, DimInfo ] \textbf{end if}
		\Statex	\hskip12pt When $\Sys$ is an {\ODE} also calculate $\mbox{LGMLin} := {\mbox{LGMLinTest}}(\Sys)$.
		
		\State Set $\hat{S}^\star :=  \left\{ \hat{\xi}^i= 0, \hat{\eta}^j_{\hat{u}^k}  = 0:  1 \leq i \leq n,  \; 1 \leq j,k \leq m  \right\}$ 
		\Statex \hskip12pt $M := S' \cup \; \hat{S}^\star|_{\Psi} \; \cup \; M_{\BK}(\mathcal{L}', \hat{\mathcal{L}}^\star ) \; \cup \{\mbox{Det} \mbox{Jac}(\Psi) \not = 0\}  $ 	
		\State Compute list of consistent cases $P = [P_1, \cdots , P_{\mbox{nc}}]$ with dim $\geq d$:
		\Statex	\hskip12pt $P := \text{\dec}(M, \prec, casesplit,  mindim=d)$ 
		\State	\textbf{if} $P=\emptyset$  \textbf{then} IsLinearizable := false \textbf{return} [ IsLinearizable, DimInfo ] \textbf{end if}	                                  
		\State		 $Q$ := [ ]
		\State		\textbf{for} $k = 1$ {\textbf{to}} $\text{nc}$ \textbf{do} 
		\State		\hskip12pt  \textbf{if} $\mbox{HF}(\Sys)  = \mbox{HF}(P_k)$ \textbf{then} $ Q := Q \cup P_k$   \textbf{end if}
		\Statex		\textbf{end do}  %
		\State		\textbf{if} $Q=$ [ ]  \textbf{then} IsLinearizable := false \textbf{return} [ IsLinearizable, DimInfo ] 
		\Statex		\hskip12pt \textbf{else if} $Q \not =$  [ ]  \textbf{then} IsLinearizable := true
		\Statex		\textbf{end if}
		\State		\textbf{if} CaseSelect $\not \in$ Options \textbf{then} $C := [ 1 ]$ 
						\textbf{else} Assign $C$ using Options 
					\textbf{end if}
		\State		\textbf{for} $c \in C$ \textbf{do} 
						$\Syshat_c := \text{ExtractTarget}(Q_c)$  
					\textbf{end do}
					(See \ref{Alg:ExtractTarget})
		\State		\textbf{for} $c \in C$ \textbf{do} 
						Attempt heuristic integration $ \Psi^c_{\mbox{sol}} :=  \maple{PDSolve}(Q_c)$   
					\textbf{end do}
					(See \ref{Method:PDSolve})

		\State		\textbf{return} $\bigcup_{c \in C}$ [$Q_c$, $\Syshat_c$, $\Syshat_c|_{{\Psi}^c_{\mbox{sol}}}$, $\Psi^c_{\mbox{sol}}$]  
		\Statex		$\hrulefill$ 
		\Statex		\textbf{Abbreviations used above}:  \dec:  Differential Elimination Completion, $M_{\BK}$ BK system,  ID: InitialData
	\end{algorithmic}
	\label{Alg:MapDE2Linear}
\end{algorithm}

\subsection{Notes on the MapDE Algorithm with Target = LinearDE}

\label{sec:NotesMapDE2Linear}

We briefly list some main aspects of Algorithm \ref{Alg:MapDE2Linear}.

\begin{itemize}
	\item[Input:]  
	Due to current limitations of {\Maple}'s \maple{DeterminingPDE} we restrict to input a single system $\Sys$, in {\dec} (i.e. {\rif}-form) leading linear equations with leading derivatives of differential order $\geq 1$, together with inequations and no leading nonlinear equations.
	This form is more general than Cauchy–Kowalevski form, and includes over and under-determined systems, but not systems with
	$0$ order (algebraic) constraints. 
	
	\maple{Options} refers to additional options for strategies and outputs.  For example including OutputDetails in 
	Options yields more detailed outputs. 
	
	\item[Step 2:] See Remark \ref{rem:DerivedAlgebraAlgorithm}, where we briefly describe our new 
	algorithm for computing determining systems for infinite dimensional derived algebras.
	\item[Step 3:]  As discussed in \S \ref{sec:PreMapEqs}, $\Sys$ being linearizable means that linear superposition 
	generates a symmetry sub-algebra of $\mathcal{L}$, yielding {\todo{[R4C11]}} some fairly well-known efficient tests 
	for rejecting many non-linearizable systems.  See Algorithm \ref{Alg:PreEquivTest} for
	details.  
	We also apply Algorithm \ref{alg:LGMLinTest} for the LGMLinTest \cite{LGM101:LG} when $\Sys$ is an {\ODE}
	in order to compare and test our Hilbert linearization test which occurs later in the algorithm.
	{\todo{[R4C12]}}
	
	\item[Step 4:]  $\hat{S}^\star|_{\Psi}$ is $\hat{S}^\star$ evaluated in $(x,u)$ coordinates via $\Psi$ using differential reduction.
	$\hat{S}^\star :=  \left\{ \hat{\xi}^i= 0, \hat{\eta}^j_{\hat{u}^k}  = 0:  1 \leq i \leq n,  \; 1 \leq j,k \leq m  \right\}$ 
	
	\item[Step 5:]  Here $\maple{mindim}= \dim(\Sys) = d$, as computed by the \Maple command \maple{initialdata} of the \dec form of $\Sys$.  The mindim option avoids computing cases of dimension $< d$
	by monitoring an upper bound based on initial data of such cases.  The block elimination ranking $\prec$ ranks all infinitesimals and their derivatives for the first block $[\xi, \eta, \hat{\xi}, \hat{\eta}]$ strictly greater 
	than the second block of the $\phi$ map variables, which are strictly greater than all derivatives of the third block of $\psi$ variables.  This maintains linearity in the variables $[\xi, \eta, \hat{\xi}, \hat{\eta}]$. The mindim dimension is computed with respect to these variables, and not the degrees of freedom in the map variables $(\psi, \phi)$.  The block structure also facilitates the later integration phase.  Each case $P_k$ consists of equations and inequations.
	
	\item[Step 9:]  The Hilbert Functions of $R$ and $P_k$, disregarding the equations that don't involve
	infinitesimals, should be equal if the system is linearizable.
	
	\item[Step 11:]  {\todo{[R1C11]}}
	Note that $Q$ can consist of several systems.  If \textit{CaseSelect = all} is included in Options, then all cases leading to linearization are returned.  By default, {\MapDE} returns only one such case: $C = [1]$.

\end{itemize}

\subsection{Proof of correctness of the MapDE Algorithm}
\label{sec:ProofMapDE}

\begin{theorem}
	\label{thm:MapDEThm}
	Let $\Sys$ be a single input system in {\rif}-form) consisting of  leading linear equations with leaders of differential order $\geq 1$, inequations, and with no leading nonlinear equations.
	Then Algorithm \ref{Alg:MapDE2Linear} converges in finitely many steps, and determines whether there exists {\todo{[R4C13]}} a local holomorphic diffeomorphism $\hat{x} = \psi(x,u), \hat{u} = \phi(x,u)$ transforming $\Sys$ to a linear homogeneous target system
	\begin{equation}
		\label{eq:Hsys}
		\Syshat: \mathcal{H} \hat{u}(\hat{x}) = 0
	\end{equation}
In the case of existence the output {\rif}-form consists of {\DPS} of equations and inequations including those for the mapping function $(\psi, \phi )$.
\end{theorem}.
\begin{proof}
We first note that Algorithm \ref{Alg:MapDE2Linear} converges in finitely many steps due to finiteness of each of the sub-algorithms used \cite{Rus99:Exi,LGM101:LG,Bou95:Rep,Robertz106:Tdec}.
To complete the proof we need to establish correctness of the two possible outcomes:  \\
\textit{Case I:} IsLinearizable = true and \textit{Case II:} IsLinearizable = false.

\vspace{0.2cm}

\noindent
\textit{\textbf{Case I: IsLinearizable = true}}  \\
Our task here is to show that given consistent input $\Sys$, IsLinearizable = true and output $Q$ then there exists a local holomorphc diffeomorphism $\Psi$ to some linear system $\Syshat$.
To do this we build initial data for a solution of $\Sys$, and initial data for solutions 
of $Q$.  A complication is that these spaces have different independent variables
$x$ and $(x,u)$.  The assumption that all leaders for the {\rif}-form of $R$ are of order $\geq 1$ enables us
to regard $(x,u)$ as independent variables for $Q$.  The inequations for $Q$ include those for $\Sys$
together with the invertibility condition $\mbox{DetJac}(\Psi) \not =0$.
	
Suppose that the input {\rif}-form of $R$ has differential order $d_R$ and consists of equations and inequations with associated varieties $ V^{=}(R)$, $ V^{\not=}(R)$ in Jet space 
$J^{d_R} (\mathbb{C}^m, \mathbb{C}^{n})$. 
So any point on the jet locus satisfies
	$(x, u^{(\leq d_R)}) \in  V^=(R) \setminus  V^{\not =}(R)$ in jet space $J^{d_R} (\mathbb{C}^m, \mathbb{C}^{n})$
where $u^{(\leq d_R)}$ denotes the jet variables of total differential order $\leq d_R$. {\todo{[R4C15]}}
Let $\pi_0^{d_R}: J^{d_R} \rightarrow X $ be the projection of points in $J^{d_R}$, the jet space of order $d_R$,
to the base space of independent variables $X \simeq \mathbb{C}^m$ where $x \in X$.
The assumption that all leaders for $R$ are of order $\geq 1$ implies that 
$\pi_0^{d_R} (V^=(R) \setminus V^{\not =} (R) )
= \mathbb{C}^m \setminus  \pi_0^{d_R} (V^{\not =} (R)) $.

When IsLinearizable = true there will be several systems $P_k$ in the list of systems $Q$ at Step 11 of Algorithm \ref{Alg:MapDE2Linear}.  We consider the case where there is only one such system, and without loss denote it by $Q$. 
For the case of several systems in $Q$ we simply repeat the argument below for each such system.
Suppose the system has differential order $d_Q$,  and consists of equations and inequations  for $v = (\xi, \eta, \psi, \phi, \hat{\eta})$ with associated varieties $V^=(Q)$, $V^{\not=}(Q)$ in $J^{d_Q}(\mathbb{C}^{m+n}, \mathbb{C}^{(2m+3n)})$, so that $((x,u), v^{(\leq d_Q)}) \in  V^=(Q) \setminus  V^{\not =}(Q)$
	in $J^{d_Q}(\mathbb{C}^{m+n}, \mathbb{C}^{(2m+3n)})$.
Then $ \pi_0^{d_Q} (V^=(Q) \setminus V^{\not =} (Q) )   = \mathbb{C}^{m+n} \setminus  \pi_0^{d_Q} (V^{\not =} (Q) )$.

Consider points $x_0 \in  \mathbb{C}^m \setminus  \pi_0^{d_R} (V^{\not =} (R)) $
and $(x_0, u_0) \in \mathbb{C}^{m+n} \setminus  \pi_0^{d_Q} (V^{\not =} (Q) )$
belonging to the projections of $R$ and $Q$ onto their base spaces $X \simeq \mathbb{C}^m$ and $X \times U \simeq \mathbb{C}^{m+n}$.  Then a family of initial data corresponding to all local analytic solutions $u$ in a neighborhood of $x_0$ exists, and similarly for $v$.
For $R$ there exists a neighborhood $\mathcal{N}(x_0, u^{(\leq d_R)}_0) \subseteq V^=(R) \setminus  V^{\not =}(R)$ in $J^{d_R}$
	and from $Q$ there exists a neighborhood  $\mathcal{N}((x_0,u_0), v^{(\leq d_Q)}_0) \subseteq  V^=(Q) \setminus  V^{\not =}(Q)$ in $ J^{d_Q}$.  The existence and uniqueness Theorems associated with {\rif}-form implies that for such analytic initial data there corresponds unique local analytic solutions and implies that there exists a local holomorphic diffeomorphism $\Psi$ between neighborhoods mapping $\Sys$ to $\Syshat$, and similarly between neighborhoods
	mapping $Q$ to {\todo{[R4C16]}} to $\hat{Q}$.  Under this diffeomorphism the images of $\mathcal{N}(x_0, u^{(\leq d_R)}_0) $
	and $\mathcal{N}((x_0, u_0), v^{(\leq d_Q)}_0)$  are $\hat{\mathcal{N}}(\hat{x}_0, \hat{u}^{(\leq d_R)}_0) $ and
	$\hat{\mathcal{N}}((\hat{x}_0,\hat{u}_0), \hat{v}^{(\leq d_Q)}_0)$  respectively.

	To show that $\Syshat$ is linear we consider the subsystem of $Q$ for $\hat{\eta}$:
	\begin{equation}
		\label{targetLin}
		\mathcal{L}^{**} = \left\{ \sum_{\ell=1}^m \hat{\eta}^\ell (\hat{x}) \dd{\hat{u}^\ell}:
		\hat{\eta}^j_{\hat{u}^k}  = 0, {\mathcal{H}} (\hat{\eta}) = 0 \right\}
	\end{equation}
where the linear system for $\hat{\eta}$ is in {\rif}-form
and ultimately we will show that ${\mathcal{H}} (\hat{\eta})$ can be taken as $\Syshat$.
First we note the linear operator $\mathcal{H}$ cannot have any coefficients depending on $\hat{u}$.  If not, and a coefficient did depend on a particular $\hat{u}^\ell$, then differentiating (\ref{targetLin}) with respect to $\hat{u}^\ell$ would yield a relation between parametric quantities, violating  the freedom to assign values independently to these parametric quantities.  This would violate the {\rif}-form and its existence and uniqueness Theorem.  So $\mathcal{H}$ only has coefficients depending on $\hat{x}$.  As a remark we note that it is 
now easily verified that $\mathcal{L}^{**} $ generates an abelian Lie pseudogroup.
	
Exponentiating the infinitesimal symmetry (\ref{targetLin}) and applying its prolongation to a solution $ \hat{u}(\hat{x})$ in
	$\hat{\mathcal{N}}(\hat{x}_0, \hat{u}^{(\leq d_R)}_0) $ yields another solution in 
	$\hat{\mathcal{N}}(\hat{x}_0, \hat{u}^{(\leq d_R)}_0) $ given by  $ \tilde{u}(\hat{x})  = \hat{u}(\hat{x}) + \hat{\eta}$
	where ${\mathcal{H}} (\hat{\eta}) = 0$.
	We have assumed in Step 9 of Algorithm \ref{Alg:MapDE2Linear} that $\mbox{HF}(\Sys)  = \mbox{HF}(Q)$, which 
	implies that all local analytic solutions in $\hat{\mathcal{N}}(\hat{x}_0, \hat{u}^{(\leq d_R)}_0) $ are of form $ \tilde{u}(\hat{x})  = \hat{u}(\hat{x}) + \hat{\eta}$ in
$\hat{\mathcal{N}}(\hat{x}_0, \hat{u}^{(\leq d_R)}_0) $.
	Consequently by a point change $\tilde{u}(\hat{x}) \rightarrow \tilde{u}(\hat{x}) - \hat{u}(\hat{x})$, 
	$\Syshat$ is equivalent to the linear homogeneous system $  {\mathcal{H}} (\hat{\eta}(\hat{x})) = 0$.  The {\rif}-form of $Q$ includes the system for $\Psi$ that determines mappings to
$\Syshat$ given by $  {\mathcal{H}} (\hat{\eta}(\hat{x})) = 0$.

{\todo{[R4C14]}}

\vspace{0.2cm}

\noindent
\textit{\textbf{Case II: IsLinearizable = false}}  \\
Suppose to the contrary that Algorithm \ref{Alg:MapDE2Linear} returns  IsLinearizable = false, yet a local analytic linearization exists.  Since the tests in Algorithm \maple{PreEquivTest} in Step 3 of Algorithm \ref{Alg:MapDE2Linear} are all necessary conditions for a linearization to exist, they all test true.
	
	Step 5 of  Algorithm \ref{Alg:MapDE2Linear} applies {\rif} using binary splitting, partitioning the jet locus into disjoint cases; and the linearization must belong to some of these cases.  By assumption, and the discussion above, there is diffeomorphism $\Psi$ of $\Sys$ in some neighborhood $\mathcal{N}(x_0, u^{(\leq d_R)}_0)$ in $J^{d_R}$ to a linear system $\Syshat$.
Further the equations of the cases corresponding to linearization in terms of
	$v = (\xi, \eta, \psi, \phi, \hat{\eta})$ must have dimension $d$.

	It cannot belong to a  case of dimension $< d$, the ones discarded by the $\mbox{mindim} = d$ option.
	Therefore by disjointness it must belong to one of the $\mbox{nc}$ cases in $P$, say $P_s$.
	Therefore this case must fail the condition that $\mbox{HF}(\Sys) = \mbox{HF}(P_s)$
	which is contrary to our assumption that such a linearization exist, completing our
	proof of correctness.
	
\end{proof}

\section{Examples}
\label{sec:Examples}
To illustrate the \MapDE Algorithm \ref{Alg:MapDE2Linear} we consider some examples.

\begin{example}(Continuation and conclusion for Examples \ref{ex:ODE3(a)},\ref{ex:ODE3(b)} and \ref{ex:ODE3(c)} using Algorithm \ref{Alg:MapDE2Linear}.)
	\label{ex:ODE3(d)}
	
	The input is (\ref{ODE3}) which is in {\rif}-form with respect to the orderly ranking 
	$u \prec u_x \prec u_{xx} \prec \cdots $, together with the inequations
	$u \not = 0, uu_{{x}}+x \not = 0,  {u}^{2}+{x}^{2} \not = 0$ or equivalently
	$u (uu_{{x}}+x)({u}^{2}+{x}^{2}) \not = 0$.
	This can be regarded as being derived from the leading linear {\DPS} which results {\todo{[R4C17]}} from 
	multiplication by factors in its denominators.
	
	\textbf{Step 1:}   Set IsLinearizable := null. Here ${\mbox{ID}}(\Sys) = [u(x_0)= c_1, u_x(x_0)= c_2,u_{xx}(x_0)= c_3]$ and $\dim \Sys = 3$, subject {\todo{[R4C18]}} to $ u (uu_{{x}}+x)({u}^{2}+{x}^{2}) \not = 0$.
	
	\textbf{Step 2:}  See Example \ref{ex:ODE3(a)} for  $S:= \mbox{\rif}(\mbox{DetSys}(\mathcal{L}))$ in (\ref{eq:ODE3DetSys}), 
	together with its $\mbox{ID}(S)$ and $\dim \mathcal{L} = \dim(S) = 4$.
	See Example \ref{ex:ODE3(c)} and in particular (\ref{eq:ODE3_L'}) for $S' := \mbox{\rif}(\mbox{DetSys}(\mathcal{L}'))$
	which yields $\dim \mathcal{L}' = \dim(S') = 3$. 
	
	\textbf{Step 3:}  Since $\dim{S} = 4 \geq d+1 = 4$ and $\dim{S'} = 3 \geq d=3$, the simplest necessary conditions for linearizability hold. 
	Also $\textbf{d} (S)  =  \textbf{d}(S')  =  \textbf{d}(R) = 0$.
	Application of Algorithm \ref{alg:LGMLinTest} for the LGMLinTest in Example \ref{ex:ODE3(b)} shows that $\Sys$ is linearizable, subject to $ u (uu_{{x}}+x)({u}^{2}+{x}^{2}) \not = 0$.
	
	
	\textbf{Step 4:}   $\mbox{Det} \mbox{Jac}(\Psi) =  \psi_x \phi_u - \psi_u \phi_x  \not = 0$, $\hat{S}^\star :=  \left\{ \hat{\xi}= 0, \hat{\eta}_{\hat{u}}  = 0  \right\}$ and $\hat{\mathcal{L}}^* $ is replaced with $\hat{\mathcal{L}}^\star$ in (\ref{eq:BKODE3}) to yield:
	\begin{equation}
		\label{eq:BKODE3repeat}
		\hspace{-1.0cm} M_{\BK}(\mathcal{L}', \hat{\mathcal{L}}^\star ) = \left\{ 
		\hat{\xi} (\hat{x}, \hat{u}) =  0  = \xi  \psi_x +  \eta \psi_u,
		\hat{\eta} (\hat{x}, \hat{u})  =   \xi  \phi_x  +\eta \phi_u  \right\}
	\end{equation}
	Evaluate $ \hat{S}^\star$ modulo $\Psi: \hat{x} = \psi(x,u), \hat{u} = \phi(x,u)$ to obtain $ \hat{S}^\star|_{\Psi}$.
	This yields $ \hat{S}^\star|_{\Psi} = \left\{ \hat{\xi} = 0, \psi_u \hat{\eta}_x - \psi_x \hat{\eta}_u = 0 \right\} $.
	Note that for brevity of notation we have replaced $\hat{\xi}(\hat{x},\hat{u})$ with $\hat{\xi}(x,u)$ and 
	$\hat{\eta}(\hat{x},\hat{u})$ with $\hat{\eta}(x,u)$.
	Thus, the mapping system $M =  S' \cup \; \hat{S}^\star|_{\Psi} \; \cup \; M_{\BK}(\mathcal{L}', \hat{\mathcal{L}}^\star ) \; \cup \{\mbox{Det} \mbox{Jac}(\Psi) \not = 0\}  $ is:
	\begin{align}
		\label{eq:Msys'}
		M = [\xi &=-{\frac {\eta\,u}{x}},\quad \eta_{{x}}={\frac { \left( {u}^{2}+{x}^{2}
				\right) \eta}{x{u}^{2}}}+{\frac {\eta_{{u}}x}{u}},   \nonumber   \\
		& \eta_{{u,u,u}}=-{
			\frac { \left( 8\,{u}^{8}+8\,{u}^{6}{x}^{2}+8\,{u}^{6}+3\,{u}^{2}+3\,{x}^{2} \right) \eta}{ \left( {u}^{2}+{x}^{2} \right) {u}^{3}}}+3\,{
			\frac {\eta_{{u}}}{{u}^{2}}},    \\
		& \hat{\xi} = 0,\quad  \psi_u \hat{\eta}_x - \psi_x \hat{\eta}_u = 0,   \quad
		\hat{\xi} =  \xi  \psi_x +  \eta \psi_u, \quad \nonumber  \\
		& \hat{\eta}  =   \xi  \phi_x  +\eta \phi_u , \quad
		\psi_x \phi_u - \psi_u \phi_x  \not = 0   ] \nonumber
	\end{align}
	
	
	\textbf{Step 5:}   Compute $P := \text{\rif}(M, \prec, casesplit,  mindim=d)$ where $d = 3$.
	This results in $3$ cases, two of which are rejected before their complete calculation since an upper bound in the computation
	drops below $mindim = d = 3$.
	The output for the single consistent case $P_1$ found is:
	\begin{align}
		\label{eq:rifMsys'}
		P_1 = [\xi &=-{\frac {\eta\,u}{x}},\eta_{{x}}={\frac { \left( {u}^{2}+{x}^{2}
				\right) \eta}{x{u}^{2}}}+{\frac {\eta_{{u}}x}{u}},   \nonumber   \\
		& \eta_{{u,u,u}}=-{
			\frac { \left( 8\,{u}^{8}+8\,{u}^{6}{x}^{2}+8\,{u}^{6}+3\,{u}^{2}+3\,{
					x}^{2} \right) \eta}{ \left( {u}^{2}+{x}^{2} \right) {u}^{3}}}+3\,{
			\frac {\eta_{{u}}}{{u}^{2}}},   \nonumber \\
		& \phi_{{x,x}}={\frac {2\,\phi_{{x,u}}x{u}^{2}-\phi_{{u,u}}{x}^{2}u+\phi
				_{{u}}{u}^{2}+\phi_{{u}}{x}^{2}}{{u}^{3}}}         , \;
		\psi_{{x}}={\frac {\psi_{{u}}x}{u}} , \\
		& \hat{\eta} =-{\frac { \left(u\phi_{{x}}-x\phi_{{u}
				} \right) \eta}{x}} , \;
		\hat{\xi}  = 0, \;   x\phi_{{u}}-u\phi_{{x}}  \neq  0 , \; \psi_u \neq  0 ]   \nonumber
	\end{align}
	
	
	\textbf{Step 6:}   $P\neq \emptyset$ contains $1$ case.
	
	\textbf{Step 7:}   Initialize $Q := [ \; ] $
	
	\textbf{Step 8:}   $k = \text{nc} = 1$
	
	
	\textbf{Step 9:}  
	{\todo{[R4C19]}}Also $\mbox{HF}(\Sys) = 1 + s + s^2$ and the {\ID} for $P_1$ yields
	$ \mbox{HF}(P_1 ) = 1 + s + s^2$.  So  $\mbox{HF}(\Sys)  = \mbox{HF}(P_1)$
	and the system is linearizable. 
	
	
	\textbf{Step 12:}  We set $Q := [P_1 ]$.  To extract the target we apply Algorithm \ref{Alg:ExtractTarget} so that 
			$\Syshat_1 := \text{ExtractTarget}(Q_1)$ which yields $\Syshat_1$ in the form:
	\begin{equation}
		\label{TargetODE3a}
		\left(\dd{\hat{x}}\right)^3 \hat{u} (\hat{x}) = a_2 (\hat{x}) \left( \dd{\hat{x}}\right)^2 \hat{u} (\hat{x}) 
		+ a_1 (\hat{x}) \dd{\hat{x}} \hat{u} (\hat{x})
		+ a_0(\hat{x})  \hat{u} (\hat{x})						
	\end{equation}
	where $a_2 (\hat{x})$, $a_1 (\hat{x})$, $a_0 (\hat{x})$ are explicit expressions in $(\hat{x},\hat{u},\hat{\psi}(\hat{x},\hat{u}), \hat{\phi}(\hat{x},\hat{u}))$ and derivatives of $\hat{\psi}(\hat{x},\hat{u}), \hat{\phi}(\hat{x},\hat{u})$.
	
	
	\textbf{Step 13:}  The $\psi$ system here is $\psi_{{x}}={\frac {\psi_{{u}}x}{u}} $.
Using $\maple{Invariants}$ from the {\LAVF} package {\todo{[R4C20]}} yields  a single invariant $x^2 + u^2$ and so $\psi = x^2 + u^2$.  Here and elsewhere the $\maple{Invariants}$ removes the need for us to specify arbitrary functions which would be the case if we started from the general solution 
of the $\psi$ equation which is in this case $\psi = F(x^2+u^2)$.
	Then substitution and solution of the $\phi$ equation then yields $\phi = G(x^2 + u^2) x + H(x^2 + u^2)$.
	The program specializes the arbitrary functions and constants to satisfy the inequations including the Jacobian condition, and in this case yields
	\begin{align}
		\label{phipsiSol1}
		\hat{x} & = \psi = x^2 + u^2    \\
		\hat{u} & = \phi = x   \nonumber 
	\end{align}
	Substitution of (\ref{phipsiSol1}) into the target (\ref{TargetODE3a}) requires first inverting (\ref{phipsiSol1}) using {\Maple}'s \maple{solve} to give $x = \hat{\psi} = \hat{u}, u = \hat{\phi} = \left( \hat{x} - \hat{u}^2 \right)^{1/2}$
	yields it explicitly as:
	\begin{equation}
		\label{TargetODE31}
		\left(\dd{\hat{x}}\right)^3 \hat{u} (\hat{x}) =
		- \frac{ (\hat{x}+1)}{\hat{x}} \hat{u} (\hat{x})						
	\end{equation}
So far our work on the integration of the mapping equations to determine the transformations is preliminary and experimental.  We have shown that the basic structure of the linear target can be determined implicitly.  It remains to
be seen how useful this would be in applications, where the mapping cannot be determined explicitly.  
Heuristic methods appear to be useful here, and we encourage the reader to try explore their own approaches.

From the output we also subsequently explored how far we could make the Target explicit before the integration of the 
map equations.  In particular we exploited the transformation (as do \cite{LGM101:LG}) that any such {\ODE} is point equivalent to one with its highest coefficients (here $a_2$, $a_1$) being zero.
	This yields additional equations on $\psi$, $\phi$ and the target takes the very simple form:
	\begin{equation}
		\label{TargetODE3b}
		\left(\dd{\hat{x}}\right)^3 \hat{u} (\hat{x})  =
		- \frac{8u^3 (u^2+x^2+1)}{(u^2+x^2) \psi_u^3} \hat{u} (\hat{x})							
	\end{equation}
	The \rif-form of the system for $\phi, \psi$ is:
	\begin{align}
		\label{rifpsiphisys}
		\psi_{{x}}&={\frac {\psi_{{u}}x}{u}}   \nonumber \\
		\psi_{{u,u,u}} &=-1/2\,{\frac {-3\,{\psi_{{u,u}}}^{2}{u}^{2}+3\,{\psi_{
						{u}}}^{2}}{\psi_{{u}}{u}^{2}}}  \nonumber \\
		\phi_{{x,x}}&=2\,{\frac {x\psi_{{u,u}} \left( \phi_{{x}}u-\phi_{{u}}x
				\right) }{\psi_{{u}}{u}^{2}}}+{\frac {\phi_{{u,u}}{x}^{2}u+\phi_{{u}}
				{u}^{2}-2\,\phi_{{x}}ux+\phi_{{u}}{x}^{2}}{{u}^{3}}}  \nonumber \\
		\phi_{{x,u}}&={\frac {
				\psi_{{u,u}}\phi_{{x}}u-\psi_{{u,u}}\phi_{{u}}x+\psi_{{u}}\phi_{{u,u}}
				x-\psi_{{u}}\phi_{{x}}}{\psi_{{u}}u}}    \nonumber  \\ 
	\end{align}
	The general solution of the system is found by {\Maple} and yields the same particular solution as before
	for $\psi, \phi$. It seems to have a made a straightforward problem, more difficult! 
	
\end{example}

{\todo{[R4C21]}}
\begin{example}(Lyakhov, Gerdt and Michels Test Set)
	Lyakhov, Gerdt and Michels \cite{LGM101:LG} consider the following 
	test set of \ODE of order $d$, for $3 \leq d \leq 15$:
	\begin{equation}
		\label{LGMTestODE}
		\left( \frac{d}{dx}\right)^d (u(x)^{2})+ u(x)^{2}=0
	\end{equation}
	By inspection this has the linearization for any $d$:
	\begin{equation}
		\label{LGMTestODEtrans}
		\Psi=\{\hat{x} = x, \, \hat{u} = u^2  \}
	\end{equation}
	{\todo{[R1C12]}}
	
They report times on an Intel(R)Xeon(R) X5680 CPU clocked at 3.33 GHz and 48GB RAM. All the following times are measured from the entry to the program.  The times for detecting the existence of the linearization by the LGM Test in \cite{LGM101:LG} range from 0.2 secs for $d = 3$ to about 150 secs for $d = 15$. For comparison, we run MapDE with our own implementation of the LGMLinTest using {\LAVF}. Our runs of the same tests to detect the existence of the linearization using on a 2.61 GHz I7-6600U processor with 16 GB of RAM range from 0.422 secs when $d = 3$ to 11.5 secs when $d= 15$.

Their method from start to existence and then construction of the linearization (existence and construction), takes 7512.9 secs for $d=9$ and is out of memory for $d \geq 10$.  In contrast, we report times for existence and construction that are only slightly longer than our existence times for $3 \leq d \leq 15$.  For $d=3$ to $d=15$ we also report the time for our Hilbert test for existence of linearization and the total time for the existence and construction of the explicit linearization. Thus LGMTest time $<$ Hilbert Test time $<$ Existence and Construction. These results are displayed in Fig.~\ref{fig:testTime} on a $\log_{10} $ axis.
	
	 \begin{figure}[h]
		\begin{centering}
			\includegraphics[angle=0,origin=c, width=6cm]{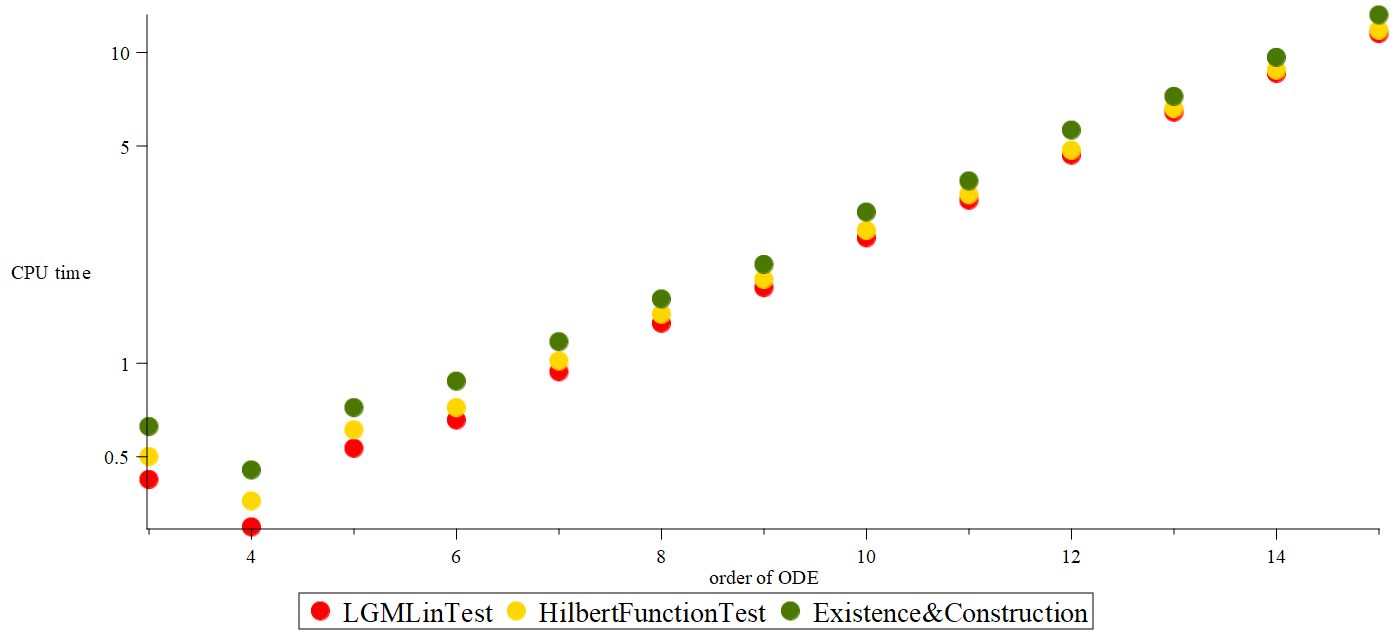}
			\caption{The graph represents the CPU times for $\left( \frac{d}{dx}\right)^d (u(x)^{2})+ u(x)^{2}=0$. Timings from $t = 0$ (start of {\MapDE}) to the time to LGMLin linearization Existence confirmation (Red), time to Hilbert Existence confirmation (Yellow), time from $t = 0$ to existence and construction of the linearizing transformations (Green).}
			\label{fig:testTime}
		\end{centering}
	\end{figure}

On this test, our approach appears to have more favorable memory behavior, which possibly is due to our equations being less nonlinear those of Lyakhov et al. \cite{LGM101:LG}. However, more testing and analysis are needed to make a reasonable comparison.    
\end{example}

\begin{example}
	\label{ex:Burgers}
	Consider Burger's equation, modeling the simplest nonlinear combination of convection and diffusion:
	\[		u_{{x,x}}= u_{{t}} -u u_x
	\]	
	Using our algorithm \MapDE with $\maple{TargetClass} = \maple{LinearDE}$ shows that it has finite dimensional
	Lie symmetry algebra with $\dim \mathcal{L} = 5 < \infty$.  Thus by the preliminary equivalence test \maple{PreEquivTest}, it is not linearizable by point transformation.  However rewriting this equation in conserved form $\dd{x} (u_x + \frac{1}{2} u^2 ) = \dd{t} u$ implies that there exists $v$: 
	\begin{equation}
		\label{potBurgers}
		v_x = u, \; \; \; v_t = u_x + \frac{1}{2} u^2 
	\end{equation}
	Applying \MapDE with $\maple{TargetClass} = \maple{LinearDE}$ to (\ref{potBurgers}) shows that this new system is linearizable, with the {\rif}-form of the $\Psi$ system given by:
	\begin{align}
		\label{mapsysBurgers}
		\phi_{{u,u}}&=0,\; \varphi_{{u,u}}=0, \; \phi_{{u,v}}=-1/2\,\phi_{{u}},  \nonumber  \\
		\varphi_{{u,v}}&=-1/2\,\varphi_{{u}}, \; \phi_{{v,v}}=-1/2\,\phi_{{v}},   \nonumber  \\
		\varphi_{{v,v}}&=-1/2\,\varphi_{{v}},
		\Upsilon_{{u}}=0,\psi_{{u}}=0,  
		\Upsilon_{{v}}=0,\psi_{{v}}=0
	\end{align}
	After integration this yields
	$\Psi$:
	\[
	{\it \xh}= \psi = x ,\, {\it \thn}= \Upsilon =  t,\, {\it \uh}= \varphi = u \exp\left(-\frac{v}{2} \right), \, {\it \vh}= \phi = \exp\left(-\frac{v}{2}\right)
	\]
	and the Target system $\Syshat$ is 
	$[{\it \uh}_{{{\it \xh}}}=-2 {\it \vh}_{{{\it \hat{t}}}},   {\it \vh}_{{{\it \xh}}}=-{\it \uh}/2]$
	so that ${\it \vh}_{\it \hat{t}} = {\it \vh}_{{{\it \xh \xh}}}$, ${\it \uh}_{\it \hat{t}} = {\it \uh}_{{{\it \xh \xh}}}$.
	This implies that the original Burger's equation is apparently also linearizable, through the introduction 
	of the auxiliary nonlocal variable $v$.  This paradox is resolved in that the resulting very useful transformation 
	is not a point transformation, since it effectively involves an integral.
	For extensive developments regarding such nonlocally related systems see \cite{Blu10:App}. 
\end{example}

\begin{example}(Nonlinearizable examples with infinite groups)
	\label{ex:KPLiouville}
	Consider the KP equation
	\begin{equation}
		\label{KP}
		u_{{x,x,x,x}}=-6\,uu_{{x,x}}-6\,{u_{{x}}}^{2}-4\,u_{{x,t}}-3\,u_{{y,y
		}}
	\end{equation}
	which has 
	\begin{align}\mbox{ID}(R) = \{ & u(x_0, y, t) = F_1(y, t),\; u_x(x_0, y, t) = F_2(y, t), \nonumber  \\ 
		&u_{x,x}(x_0, y, t) = F_3(y, t),\; u_{x,x,x} (u)(x_0, y, t) = F_4(y, t)\}.
	\end{align}
	
	Applying {\MapDE} shows that the defining system $S$ for symmetries $\xi \dd{x} + \eta \dd{y} + \tau \dd{t} + \beta \dd{u}$, has initial data which is the union of infinite data along the Hyperplane $p = (x_0, y_0, t, u_0)$ and finite 
	initial data at the point $z_0 = (x_0, y_0, t_0, u_0)$:
	\begin{align}
		\mbox{ID}(S) &= \{ \beta \left( p \right) =H_{{1}} \left( t
		\right) , \beta_y \left( p \right) =H_{{2}} \left( t \right) , \beta_{y,y}  \left( p \right) =H_{{3
		}} \left( t \right)  \}     \nonumber   \\
		\phantom{XX} &   \hspace{3cm} \cup     \\
		\{ \beta_x  (z_0 ) &=c_{{1}}, \beta_u  (z_0 ) =c_{{2}},\eta (z_0 ) =c_{{3}},\eta_t  ( z_{{0}}) =c_{{4}},\tau (z_0 ) =c_{{5}},\xi (z_0 ) =c_{{6}}  \}   \nonumber
	\end{align}
	So both the KP equation and its symmetry system have infinite dimensional solution spaces:
	$\dim R = \dim S = \infty$ since both have arbitrary functions in their data.
	However the KP equation has differential dimension $\textbf{d} (R) = 2$ since there are a max of $2$ free variables ($y,t$) in its initial data while its symmetry system has  $\textbf{d} (S)  =  1$ since it has a max
	of one free variable in its data.  Thus $\textbf{d} (S)  =  1 <  2 = \textbf{d} (R)$ and by the \maple{PreEquivTest} the KP equation is not linearizable by point transformation.
	
	Consider Liouville's equation $u_{{x,x}} + u_{{y,y}} = \mbox{e}^u$ which we rewrite as a {\DPS} using {\Maple}'s function \maple{dpolyform}.  That yields $v = \mbox{e}^u$ and the Liouville equation in the form
	$v_{{x,x}}  =  - v_{{y,y}}  +  v_x^2/v + v_y^2/v + v^2$.
	{\MapDE} determines that $\dim \mathcal{L} = \infty = \dim \Sys$, and also that the Liouville equation is not linearizable by point transformation.   Interestingly it is known that Liouville's equation is linearizable by contact transformation (a more general transformation involving derivatives).
	For extensive developments regarding such contact related systems see \cite{Blu10:App}.
	
\end{example}

\begin{example}
	\label{ex:NLTelegraph}
	Given that exactly linearizable systems are not generic among the class of nonlinear systems, a natural question is how to identify such linearizable models. {\todo{[R3C2]}}Since linearizability requires large (e.g. $\infty$ dimensional) symmetry groups a natural approach is to embed a model in a large class of systems and seek the members of the class with the largest symmetry groups.
	Indeed in Wittkopf and Reid 
	\cite{ReWit113:rif} developed such as approach. We now illustrate how this approach
	can be used here. For a nonlinear telegraph equation one might embed it in the general class of spatially dependent nonlinear 
	telegraph systems
	\begin{equation}
		\label{polTel}
		v_x = u_t, \; \; \; v_{t}  = C(x,u)u_x + B(x,u)
	\end{equation}
	where $B_u \neq 0, C_u \neq 0, B_x \neq 0, C_x \neq 0$.
	Then applying the \maple{maxdimsystems} algorithm  available in {\Maple} with $dim = \infty$ a quick calculation yields
	$11$ cases, only $4$ of which satisfy the dimension restriction which we further narrow by requiring restriction to those that  have the greatest freedom in $C(x,u)$, $B(x,u)$.
	Integration yields the linearizable class:
	\begin{equation}
		\label{TelLin}
		v_x = u_t, \; \; \; v_{t}  =  \frac{1}{q_x u}  f\left( \frac{u}{q_x}\right) u_x  - \frac{q_{xx} } {q_x^2 }  f\left( \frac{u}{q_x}\right)
	\end{equation}
	and the linearizing transformation
	\begin{equation}
		\label{TelLinTran}
		\hat{x} = \frac{u}{q_x}, \hskip6pt  \hat{t} = v, \hskip6pt \hat{u} = q(x), \hskip6pt \hat{v} = t
	\end{equation}
	Similarity, we can seek the maximal dimensional symmetry group for the normalized linear Schr\"odinger Equation {\todo{[R3C3]}}
	\begin{equation}
		\label{ex:ConsCoeffR}	
		i  \hslash  \varphi_t = -\frac{\hslash^2}{2 m}
		\nabla^2  \varphi +V(x,y,t) \varphi 
	\end{equation}
	Restricting to $2$ space plus one time yields
	$V(x,y,t) = \omega(t) (x^2 + y^2) + b(t)x + c(t)y + d(t)$, and satisfies the conditions for mapping to a constant coefficient {\DE} via the methods of our previous paper \cite{MohReiHua19:Intro}.
	
\end{example}
\section{Discussion}
\label{sec:Discussion}



In this paper we give an algorithmic extension of {\MapDE} introduced in 
\cite{MohReiHua19:Intro}, that decides whether an input {\DPS} can be mapped by local holomorphic diffeomorphism to a linear system, returning equations for the mapping in {\rif}-form, useful for further
applications.
This work is based on creating algorithms that exploit results due to 
Bluman and Kumei \cite{BK105:BluKu,Blu10:App} and some aspects of \cite{LGM101:LG}.
This is a natural partner to the algorithm for deciding the existence of an invertible map of a linear {\DPS} to a constant coefficient linear {\DE} given in our previous paper \cite{MohReiHua19:Intro}.


The mapping approach \cite{LGM101:LG} for {\ODE} explicitly introduces a target linear system $\Syshat$ with undetermined coefficients, then uses the full nonlinear determining equations for the mapping and applies the \maple{ThomasDecomposition} Algorithm~\cite{Thomas106:Tdec,Robertz106:Tdec}. 
In contrast, like Bluman and Kumei, we exploit the fact that the target appears implicitly as a subalgebra of the Lie symmetry algebra $\mathcal{L}$ of $\Sys$ and avoid using the full nonlinear determining equations
for the transformations. 
Unlike Bluman and Kumei, who depend on extracting this subalgebra by explicit non-algorithmic integration, we use algorithmic differential algebra.
We exploit the fact that the subalgebra corresponding to linear super-position appears as a 
subalgebra of the derived algebra $\mathcal{L}'$ of $\mathcal{L}$, generalizing the technique for {\ODE} in \cite{LGM101:LG}.
Instead of using the {\BK} mapping equations~\eqref{eq:BK},  \cite{LGM101:LG} apply the transformations directly to the \ODE.  In contrast, our method works at the linearized Lie algebra level instead of the nonlinear Lie Group level used in  \cite{LGM101:LG}  which may be a factor in the increased space and time usage for their test set compared to our timings  given in Fig.~\ref{fig:testTime}.

Bluman and Kumei give necessary conditions in \cite[Theorem 2.4.1]{Blu10:App} 
and sufficient conditions \cite[Theorem 2.4.2]{Blu10:App} for linearization of nonlinear {\PDE} systems with $m \geq 2$.
Their requirement $\dim{\mathcal{L}} = \infty$ is dropped in our approach allowing us to deal with {\ODE} and also linearization of overdetermined {\PDE} systems.
They also use the Jacobian condition to introduce a solved form of the BK equations with coefficients 
$\alpha^i_\sigma(x,u), \beta^\nu_\sigma(x,u)$ (see \cite[Eq 2.69]{Blu10:App}), and further decompose the resulting system with respect to their $f^\sigma (\phi)$'s (our $\hat{\eta}$'s).
This decomposition results for input {\PDE} having no zero-order (i.e. algebraic) relations among the input systems of {\PDE}, a condition that is not explicitly given in the hypotheses of their theorems.   
We are planning to take advantage of this decomposition in future work, as an option to {\MapDE}, since it can improve efficiency, when applicable.
For the more general case of contact transformations for $m=1$, not considered here, see
\cite[Theorems 2.4.3-2.4.4]{Blu10:App}.

An important aspect of Theorem \ref{thm:MapDEThm} concerning the correctness of {\MapDE}, is
to show that the existence of an infinitesimal symmetry $\sum_\ell \hat{\eta}^\ell (\hat{x}) \dd{\hat{u}^\ell}$
where ${\mathcal{H}} (\hat{\eta}) = 0$, when exponentiated to act on a local analytic solution of $\Syshat$ produces all local analytic solutions in a neighborhood.
Showing this depends on showing 
$\mbox{HF}(\Sys)  = \mbox{HF}(P_k)$ in Step 9 of Algorithm \ref{Alg:MapDE2Linear} or in intuitive terms,
the size of solution space of the input system is the same as the size of the symmetry subgroup corresponding to $P_k$.
In contrast the statement and proof of \cite[Theorem 2.4.2]{Blu10:App} appears to miss this crucial hypothesis about the size of the 
solution space of the input system being equal to the size of the solution space of the symmetry sub-group
for linearization (as measured by Hilbert Series).

It is important to develop simple, efficient tests to reject the existence of mappings, based on structural and dimensional information. In addition to existing tests \cite{LGM101:LG},
\cite{Blu10:App,AnBlWo110:Mapp,TWolf108:ConLa} we introduced a refined dimension test based on Hilbert Series.  We will extend these tests in future work.
We note that the potentially expensive change of rankings needed by our algorithms (for example to determine 
the derived algebra when it is infinite dimensional) could be more efficiently accomplished by the change of rankings approach given in \cite{BouLemMazChangeOrder:2010}.


{\todo{[R4C22]}}
 Mapping problems such as those considered in this paper are theoretically and computationally challenging. Given that nonlinear systems are usually not linearizable, a fundamental problem is to identify such linearizable models. For example \cite{AnBlWo110:Mapp,TWolf108:ConLa} use multipliers for conservation laws to facilitate the determination of linearization mappings. Wolf's approach \cite{TWolf108:ConLa} enables the determination of partially linearizable systems. Setting up such problems by finding an appropriate space to define the relevant mappings is important for discovering new non-trivial mappings.
See Example \ref{ex:Burgers} and \cite{Blu10:App}
for such embedding approaches where the model is embedded in spaces that have a natural relation to the original space in terms of solutions but are not related by invertible point transformation.
Another method is to embed a given model in a class of models and then efficiently seek the members of the class with the largest symmetry groups and most freedom in the functions/parameters of the class. Example \ref{ex:NLTelegraph} illustrates this strategy. 

We provide a further integration phase to attempt to find the mappings explicitly, based on {\Maple}'s \maple{pdsolve} which will be developed in further work.{\todo{[R4C9]}}
Even if the transformations can't be determined explicitly, they can implicitly identify important features.
Linearizable systems have a rich geometry that we are only beginning to exploit, such as the availability group action on the source and target. This offers interesting opportunities to use invariantized methods, such as invariant differential operators, and also moving frames \cite{Man10:Pra,Hub09:Dif,Fel99:Mov,Arnaldsson17:InvolMovingFrames,LisleReid2006}.
Furthermore, they are available for the application of symbolic and symbolic-numeric approximation methods, a possibility that we will also explore.
Finally a model that is not exactly linearizable may be \textit{close} to a linearizable model or other attractive target, providing motivation for our future work on approximate mapping methods.


\section*{Acknowledgments}
One of us (GR) acknowledges his debt to Ian Lisle who was the initial inspiration behind this work,
especially his vision for Lie pseudogroups and an algorithmic calculus for symmetries of differential systems. GR and ZM acknowledge support from GR's NSERC discovery grant.

Our program and a demo file included some part of our computations are publicly available on GitHub at: 
\href{https://github.com/GregGitHub57/MapDETools}{https://github.com/GregGitHub57/MapDETools}





\addcontentsline{toc}{section}{Index}
\label{ind:index}
\printindex

\end{document}